\input graphicx
\input amssym.def
\input amssym
\magnification=1000
\baselineskip = 0.20truein
\lineskiplimit = 0.01truein
\lineskip = 0.01truein
\vsize = 8.7truein
\voffset = 0.1truein
\parskip = 0.10truein
\parindent = 0.3truein
\settabs 12 \columns
\hsize = 6.0truein
\hoffset = 0.1truein
\font\ninerm=cmr9

\setbox\strutbox=\hbox{%
\vrule height .708\baselineskip
depth .292\baselineskip
width 0pt}
\font\caps=cmcsc10
\font\smallcaps=cmcsc9
\font\bigtenrm=cmr10 at 14pt

\font\smallsl=cmsl9

\def\sqr#1#2{{\vcenter{\vbox{\hrule height.#2pt
\hbox{\vrule width.#2pt height#1pt \kern#1pt
\vrule width.#2pt}
\hrule height.#2pt}}}}

\centerline{\bigtenrm ELEMENTARY KNOT THEORY}
\tenrm
\vskip 14pt
\centerline{MARC LACKENBY}
\vskip 18pt

\vskip 18pt
\centerline{\hbox to 5.5truein {1. Introduction \dotfill 1}}
\centerline{\hbox to 5.5truein {2. The equivalence problem \dotfill 1}}
\centerline{\hbox to 5.5truein {3. The recognition problem \dotfill 8}}
\centerline{\hbox to 5.5truein {4. Reidemeister moves \dotfill 15}}
\centerline{\hbox to 5.5truein {5. Crossing number \dotfill 19}}
\centerline{\hbox to 5.5truein {6. Crossing changes and unknotting number \dotfill 23}}
\centerline{\hbox to 5.5truein {7. Special classes of knots and links \dotfill 25}}
\centerline{\hbox to 5.5truein {8. Epilogue \dotfill 29}}

\vskip 18pt
\centerline{\caps 1. Introduction}
\vskip 6pt

In the past 50 years, knot theory has become an extremely well-developed subject.
But there remain several notoriously intractable problems about knots and links, 
many of which are surprisingly easy to state. The focus of this article is this 
`elementary' aspect to knot theory. Our aim is to highlight what we still do not understand,
as well as to provide a brief survey of what is known. The problems that we will concentrate on are the
non-technical ones, although of course, their eventual solution is likely to be sophisticated. 
In fact, one of the attractions of knot theory is its extensive interactions with
many different branches of mathematics: 3-manifold topology, hyperbolic geometry, 
Teichm\"uller theory, gauge theory, mathematical physics, geometric group theory, graph theory
and many other fields. 

This survey and problem list is by no means exhaustive. For example, we largely neglect the theory
of braids, because there is an excellent survey by Birman and Brendle [5] on this topic. We have also avoided
4-dimensional questions, such as the slice-ribbon conjecture (Problem 1.33 in [41]). Although these do have a significant
influence on elementary knot theory, via unknotting number for example, this field is so extensive
that it would best be dealt with in a separate article.

In fact, this article comes with a number of health warnings. It is frequently vague, often deliberately so,
and many standard terms that we use are not defined here. It is certainly far from complete, and we apologise
to anyone whose work has been omitted. And it is not meant to be a historical account of the many developments
in the subject. Instead, it is primarily a list of fundamental and elementary problems,
together with enough of a survey of the field to put these problems into context.

\vskip 18pt
\centerline {\caps 2. The equivalence problem}
\vskip 6pt

The equivalence problem for knots and links asks the most fundamental question in the field: {\sl can we decide
whether two knots or links are equivalent?} Questions of this sort arise naturally in just about any
branch of mathematics. In topology, there are many well-known negative results. A central result of
Markov [65] states that there is no algorithm to determine whether two closed $n$-manifolds are homeomorphic,
when $n \geq 4$. But in dimension 3, the situation is more tractable. In particular, the equivalence
problem for knots and links is soluble. In fact, there are now five known ways to solve this problem,
which we describe briefly below.\footnote{$^1$}{\ninerm Since the first version of this article was produced,
two new solutions to the related problem of determining whether
two closed 3-manifolds are homeomorphic have been given by Scott and Short [91] and by Kuperberg [49].
Both rely on the solution to the Geometrisation Conjecture, as in Section 2.2.}
Four of these techniques focus on the {\sl exterior} $S^3 - {\rm int}(N(K))$
of the link $K$. Now, one loses a small amount information when passing to the link exterior,
because it is possible for distinct links to have homeomorphic exteriors. (This is not the case
for knots, by the famous theorem of Gordon and Luecke [23]). So, one should also keep track
of a complete set of meridians for the link, in other words, a collection of
simple closed curves on $\partial N(K)$, one on each component of $\partial N(K)$,
each of which bounds a meridian disc in $N(K)$. Then, for two links $K$ and $K'$,
there is a homeomorphism between their exteriors that preserves these meridians
if and only if $K$ and $K'$ are equivalent.

\vskip 6pt
\noindent  {\caps 2.1. Hierarchies and normal surfaces} 
\vskip 6pt

Haken [25], Hemion [30] and Matveev [66] were the first to solve the equivalence problem for links. Their
solution is lengthy and difficult. A full account of their argument was given only fairly recently by Matveev [66].

Haken used incompressible surfaces arranged into sequences called hierarchies. Recall that a surface $S$ properly
embedded in a 3-manifold $M$ is {\sl incompressible} if, for any embedded disc $D$ in $M$
with $D \cap S = \partial D$, there is a disc $D'$ in $S$ with $\partial D' = \partial D$.
Given such a surface, it is natural to cut $M$ along $S$,  creating the
compact manifold ${\rm cl}(M - N(S))$. It turns out that one can always find another
properly embedded incompressible surface in this new manifold, and thereby iterate
this procedure. A {\sl hierarchy} is a sequence of manifolds and surfaces
$$M = M_1 \buildrel S_1 \over \longrightarrow M_2 \buildrel S_2 \over \longrightarrow \dots
\buildrel S_{n-1} \over \longrightarrow M_n$$
where each $S_i$ is a compact orientable incompressible surface properly embedded
in $M_i$ with no 2-sphere components, where $M_{i+1}$ is obtained from $M_i$
by cutting along $S_i$ and where the final manifold $M_n$ is a collection of 3-balls. Haken proved that
the exterior of a non-split link always admits a hierarchy. He was able to solve the
equivalence problem for links via the use of hierarchies with certain nice properties.
An example of a hierarchy for a knot exterior is given in Figure 1.

Associated with a hierarchy, there is a 2-complex embedded in the manifold $M$.
To build this, one extends each surface $S_i$ in turn, by attaching a collar to $\partial S_i$ in
$N(S_1 \cup \dots \cup S_{i-1})$, so that the boundary of the new surface runs
over the earlier surfaces and $\partial M$. The union of these extended surfaces and $\partial M$
is the 2-complex. The 0-cells are the points where three surfaces intersect.
The 1-skeleton is the set of points which lie on at least two surfaces. Since the
complement of this extended hierarchy and the boundary of the manifold is a collection
of open balls, we actually obtain a cell structure for $M$ in this way. It is a trivial but
important observation that this is a representation of $M$ that depends only on the topology
of the hierarchy. Haken's solution to the equivalence problem establishes that, as long as one
chooses the hierarchy correctly, there is a way of building this cell structure,
starting only from an arbitrary triangulation of $M$.

\vskip 12pt
\centerline{
\includegraphics[width=3.7in]{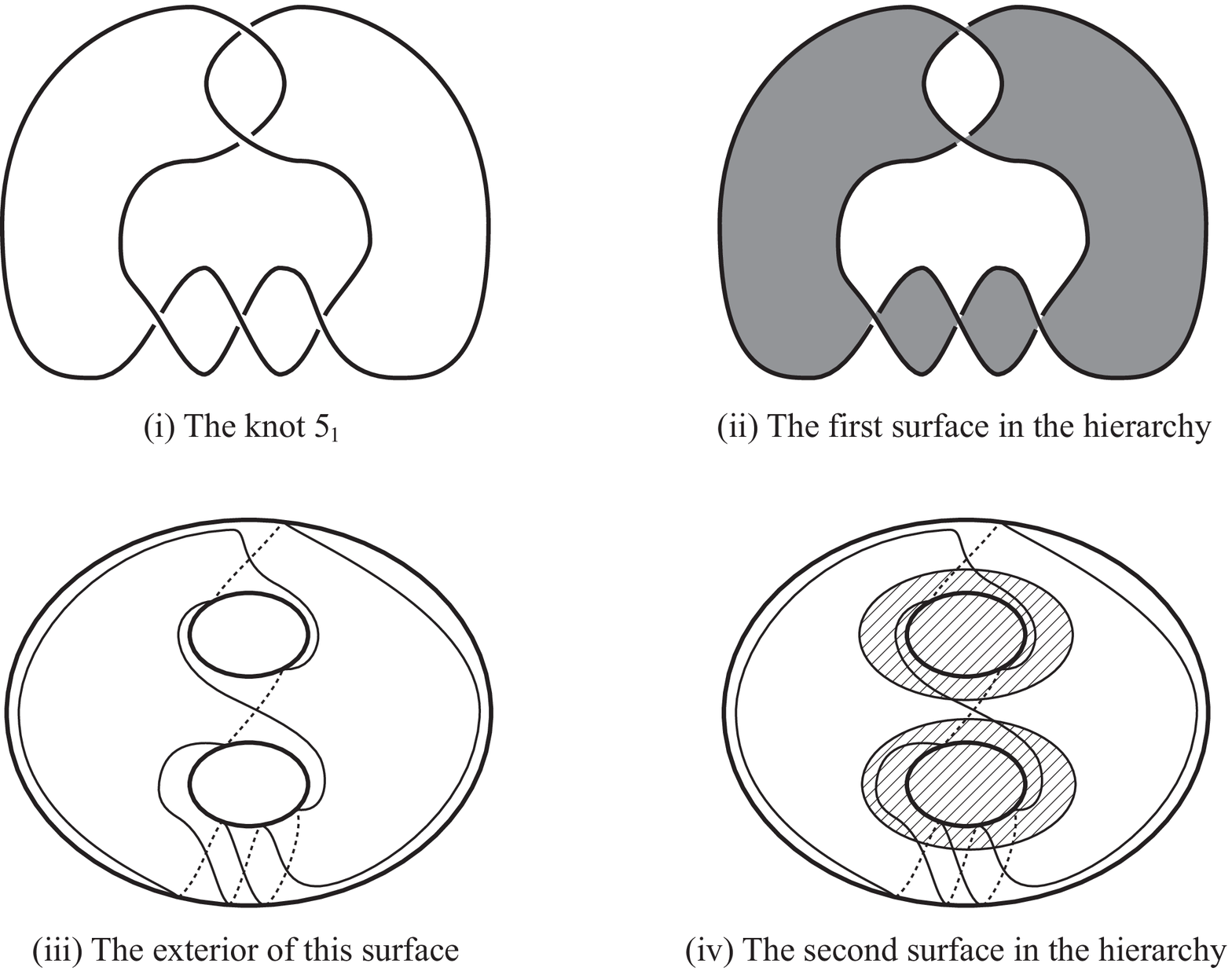}
}
\vskip 6pt
\centerline{Figure 1: A hierarchy for a knot exterior}

Suppose that we are given diagrams for two (non-split) links $K$ and $K'$. The first step in the
algorithm is to construct triangulations $T$ and $T'$ for their exteriors $M$ and $M'$.
We know that $M$ and $M'$ admit hierarchies with the required properties.
Moreover, if $M$ and $M'$ are homeomorphic (by a homeomorphism preserving
the given meridians), then this homeomorphism takes one such hierarchy for $M$
to a hierarchy for $M'$. So, let us focus initially on the hierarchy for $M$.
The next stage in the procedure is to place the first surface $S$ of the hierarchy into {\sl normal form}.
By definition, this means that it intersects each tetrahedron of $T$ in a collection of
triangles and squares, as shown in Figure 2.

\vskip 6pt
\centerline{
\includegraphics[width=3in]{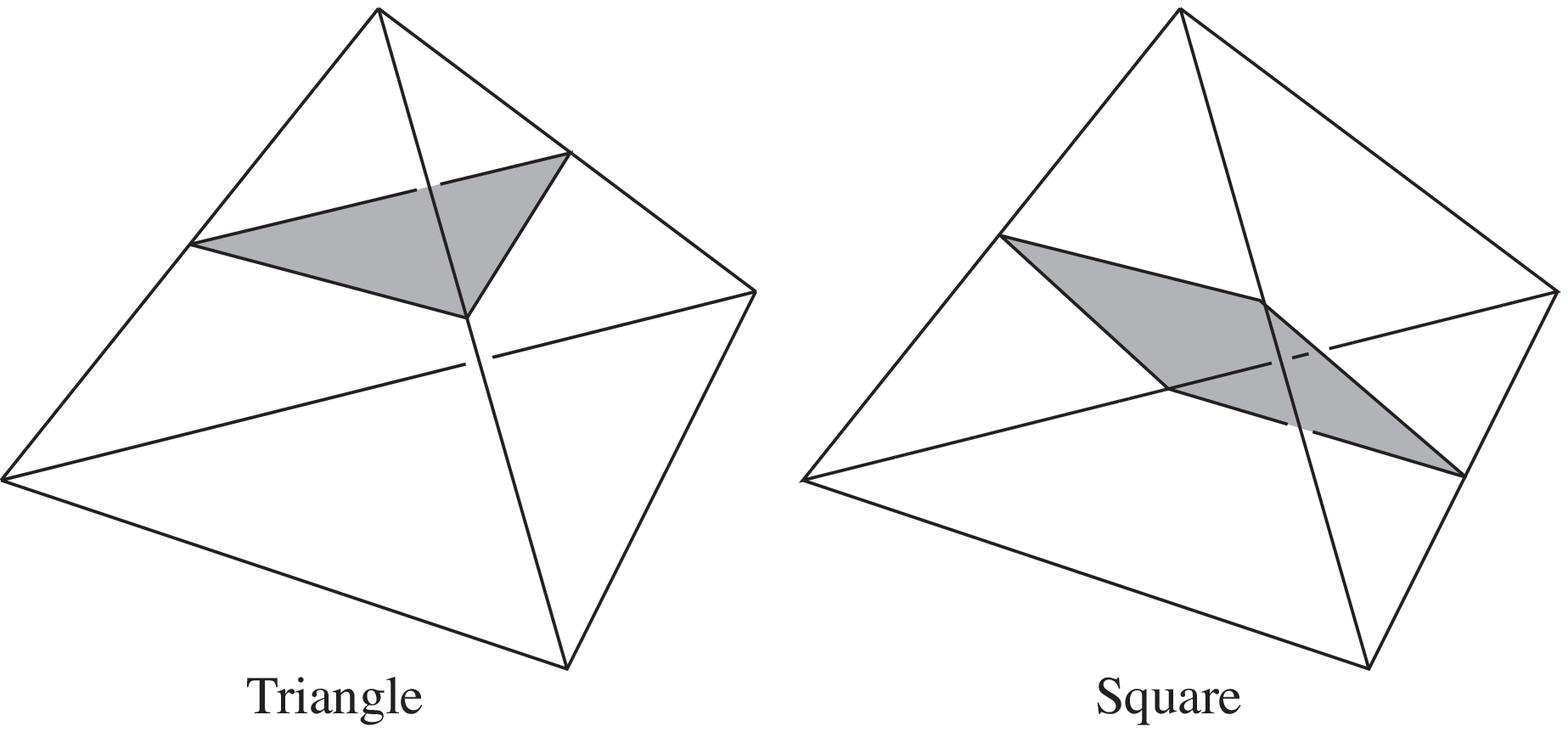}
}
\vskip 6pt
\centerline{Figure 2: Normal surface}

Now, any incompressible, boundary-incompressible surface $S$ (with no components that
are 2-spheres or boundary-parallel discs)
in a compact orientable irreducible 
boundary-irreducible 3-manifold $M$ with a given triangulation $T$ may be ambient isotoped into normal form.
It can then be encoded by a finite amount of data. The normal surface $S$ determines a collection of $7t$ integers,
where $t$ is the number of tetrahedra in $T$. These integers simply record
the number of triangles and squares of $S$ in each tetrahedron. (Note from Figure 2 that there are
$4$ types of normal triangle and $3$ types of normal square in each tetrahedron.) The list of
these integers is denoted $[S]$ and is called the associated {\sl normal surface vector}.
Normal surface vectors satisfy a number of constraints. Firstly, each co-ordinate
is, of course, a non-negative integer. Secondly, they satisfy various linear equations
called the {\sl matching conditions} which guarantee that, along each face of $T$ with a tetrahedron on each side,
the triangles and squares of $S$ in the two adjacent tetrahedra patch together correctly.
Finally, there are the {\sl quadrilateral constraints} which require $S$ to have
at most one type of square in each tetrahedron. This is because two different square
types cannot co-exist in the same tetrahedron without intersecting. Haken showed
that there is a one-one correspondence between properly embedded normal surfaces (up to a suitable
notion of `normal' isotopy) and vectors satisfying these conditions.
Given this strong connection, it is natural to use tools from linear algebra when studying normal surfaces. In
particular, one can speak of a normal surface $S$ as being the {\sl sum} of
two normal surfaces $S_1$ and $S_2$ if $[S] = [S_1] + [S_2]$. A normal surface
$S$ is called {\sl fundamental} if it cannot be expressed as a sum of two other
non-empty normal surfaces. Haken showed that there is a finite list of fundamental normal
surfaces, and these are all algorithmically constructible. A key part of Haken's argument
is to show that some surfaces can be realised as fundamental normal surfaces
with respect to {\sl any} triangulation of the manifold. Thus, if we could find
a hierarchy where each surface had this property, then we could construct the
hierarchy starting with any triangulation $T$ of $M$. In fact, Haken did not prove
that there is a hierarchy with this property, but he did show that there is a hierarchy
where each surface can be realised as a sum of a bounded number of fundamental normal
surfaces. This is sufficient to make the argument work.

Thus, Haken's algorithm essentially proceeds by starting with the triangulations $T$ and $T'$
of the two link exteriors. Then, one constructs all possible hierarchies for each
manifold, with the property that each surface in the hierarchy is a bounded sum
of fundamental surfaces. From each such hierarchy, one forms the associated
cell structure for $M$. If two such cell structures, one from each of $T$ and $T'$,
are combinatorially equivalent, then the links are the same. If none of these cell structures
are equivalent, then the links are distinct.

There were two cases that Haken could not handle using these methods.
When $M$ fibres over the circle with fibre $S$, then cutting $M$ along $S$
results in a copy of $S \times I$. In $S \times I$, there is no good choice of
surface to cut along. In particular, it is not clear how to find a surface that is guaranteed to be fundamental
or at least a bounded sum of fundamental surfaces. The resolution of the equivalence problem for these
manifolds was not completed until work of Hemion [30], where he gave a solution to the
conjugacy problem in the mapping class group of $S$. A similar problem arises
when cutting along $S$ yields a twisted $I$-bundle (in which case, $S$ is known
as a {\sl semi-fibre}). This case was fully dealt with by Matveev [66].

In fact, there is a
now an alternative way of sidestepping this issue, by establishing that there is always 
a hierarchy for a non-split link exterior which is always constructible. This relies
on the result of Culler and Shalen [13], which says that if $M$ is a compact orientable hyperbolic 3-manifold
with boundary a non-empty collection of tori, then $M$ contains either a closed, non-separating, essential
properly embedded surface or a separating, connected, essential, properly embedded surface with non-empty boundary.
When $M$ is the exterior of a hyperbolic link in the 3-sphere, this surface can be
chosen to be neither a fibre nor a semi-fibre (see Section 7 in [9]). By work of Mijatovic [71], 
there is such a surface that can be realised as a bounded sum of fundamental
surfaces.

Of course, this summary is a dramatic over-simplification. The details of the argument
are very complicated. A full exposition, which lasts several hundred pages, can be
found in the excellent book by Matveev [66].

\vskip 6pt
\noindent {\caps 2.2. Geometric structures}
\vskip 6pt

One of the most striking advances in low-dimensional topology was Thurston's introduction of
hyperbolic geometry, via his Geometrisation Conjecture [100]. This was proved by Perelman in 2003 [81, 82, 83],
but Thurston himself proved his conjecture in the case of knot and link complements,
and this can be used to provide another method of solving the equivalence problem.
An approximate statement of the Geometrisation Conjecture is that every compact orientable
3-manifold admits a `canonical decomposition' into `geometric pieces'.

The canonical decomposition takes place in two steps. (We focus for simplicity on the
case where the boundary of the 3-manifold is empty or a collection of tori.) Firstly, one decomposes
the 3-manifold into its prime connected summands. Secondly, one cuts the manifold
along its JSJ tori. Both processes can be achieved algorithmically using normal surface theory [34].
(See [75] and [66] for more details of JSJ theory, which was first developed by Jaco, Shalen [33]
and Johannson [36].)

Once this decomposition has taken place,
the resulting pieces either are Seifert fibred or admit finite-volume hyperbolic structures. The generic situation is the hyperbolic
case, and it is this that we will focus on. An algorithm to find a finite-volume hyperbolic structure on a
3-manifold, provided one exists, was given by Manning [63], based on work of Casson.
By combining [63] with work of Weeks [103], it possible to compute the Epstein-Penner
decomposition of the manifold [18]. This is a striking construction: it is a decomposition of
a non-compact finite-volume hyperbolic 3-manifold into ideal polyhedra. The crucial
thing about it is that it is canonical: it depends only on the topology of the manifold,
not on any arbitrary choices. Thus, two finite-volume non-compact hyperbolic 3-manifolds are homeomorphic
if and only if their Epstein-Penner decompositions are combinatorially equivalent.
Moreover, one can determine whether there is a homeomorphism between the
manifolds, taking one given collection of disjoint simple closed curves in the boundary
to another. This therefore gives a solution to the equivalence problem for knots and links.

This is more than just a theoretical algorithm. The computer programs Knotscape [31],
Snappea [102] and Snap [10] all attempt to compute the Epstein-Penner decomposition of
a hyperbolic 3-manifold. Unlike the Casson-Manning algorithm, they are not guaranteed
to do so. But {\sl if} Snap or Knotscape finds what it claims is the Epstein-Penner decomposition, then this is 
indeed the correct decomposition,
and as a result, these programs are a practical method of reliably determining whether two links
are equivalent.

\vskip 6pt
\noindent {\caps 2.3. Relatively hyperbolic groups}
\vskip 6pt

One can also use some of the machinery of geometric group theory to solve the
equivalence problem, at least for hyperbolic knots. 

The fundamental group of a closed hyperbolic 3-manifold is Gromov-hyperbolic. 
Mostow rigidity states that two such manifolds are isometric if and only
if their fundamental groups are isomorphic. Moreover, Sela provided an algorithm
to decide whether two torsion-free Gromov-hyperbolic groups are isomorphic [92]. 
Now, the complement of a knot is, of course, not closed.
But when this complement has a hyperbolic structure, its fundamental group is a relatively
hyperbolic group, relative to the fundamental group of the toral boundary.
Mostow rigidity applies in this case also. Sela's result was extended to this
class of relatively hyperbolic groups by Dahmani and Groves [14]. So, at least in 
the case of hyperbolic knots, one can solve the equivalence problem this way.
This does not immediately lead to a solution for hyperbolic links, because links
are not determined by their complements. But presumably, the algorithm of Dahmani
and Groves could be adapted to take account of a set of meridians in the boundary tori.
This method also only works in the hyperbolic case, and so to turn this into a fully-fledged
solution to the equivalence problem for all links, presumably one would first need to
find the `canonical decomposition into geometric pieces' that was discussed in
the previous subsection. Nevertheless, this is a genuinely different solution to
the equivalence problem. But its algorithmic complexity seems to be hard to
estimate.

\vskip 6pt
\noindent {\caps 2.4. Reidemeister moves}
\vskip 6pt

There is an alternative way of interpreting the equivalence problem in terms of Reidemeister moves. Recall
that a {\sl Reidemeister move} is a local modification to a link diagram as shown in Figure 3.
Reidemeister proved [85] that any two diagrams of a link differ by a sequence of Reidemeister moves.
The algorithmic significance of this is that {\sl if} two diagrams represent the same link, then
one may always prove this, given enough time and computing power. This is because one can
apply all possible Reidemeister moves to the first diagram, and thereby obtain a collection
of diagrams of the first link. One can then apply all possible Reidemeister moves to each of these,
and so on. It is a consequence of Reidemeister's theorem that if the two initial diagrams
represent the same link, then this process is guaranteed to produce a sequence of 
Reidemeister moves taking one diagram to the other.

\vskip 12pt
\centerline{
\includegraphics[width=4in]{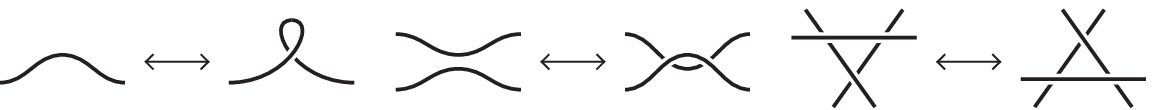}
}
\vskip 6pt
\centerline{Figure 3: Reidemeister moves}

Of course, this does not provide a solution to the equivalence problem, because if two
diagrams represent distinct link types, then the above process does not terminate. But
if one knew in advance how many Reidemeister moves were required to take one diagram
to other, then one could stop the process when one had tried all possible sequences
of Reidemeister moves of this length and if a sequence of Reidemeister moves taking one diagram
to the other had not been found, then one could declare that the links are distinct.
More specifically, a computable upper bound on the number of Reidemeister moves
required to relate two diagrams of the same link provides a solution to the equivalence
problem. In fact, it is not hard to show that the converse is also true: if there is a solution
to the equivalence problem, then, given natural numbers $n_1$ and
$n_2$, one can compute an upper bound on the number of Reidemeister
moves required to relate two diagrams of a link with $n_1$ and $n_2$ crossings.
One enumerates all link diagrams with these numbers of crossings,
and then one sorts them into their various link types, using the hypothesised
algorithm to solve the equivalence problem. Then, for all diagrams of the same
link type, one starts applying Reidemeister moves to these diagrams. By Reidemeister's
theorem, eventually a sequence of such moves will be found relating any two diagrams
of the same link. Hence, one can compute an upper bound on the number of 
moves that are required.

Recently, Coward and the author [9] have provided an explicit, computable 
upper bound on Reidemeister moves, thereby giving a new
and conceptually simple solution to the equivalence problem. Unfortunately, it is
a huge bound, and so the resulting algorithm is very inefficient. We will provide more
details in Section 4, which is devoted to Reidemeister moves.

\vskip 6pt
\noindent {\caps 2.5. Pachner moves}
\vskip 6pt

It is well known that two triangulated manifolds are PL-homeomorphic if and only if they
differ by a sequence of Pachner moves [80]. These are modifications that change the triangulation
in a very simple way. In the case of closed 3-manifolds, the moves are shown in Figure 4.
For 3-manifolds with boundary, there are some extra moves which change the triangulation
on the boundary.

\vskip 18pt
\centerline{
\includegraphics[width=4.5in]{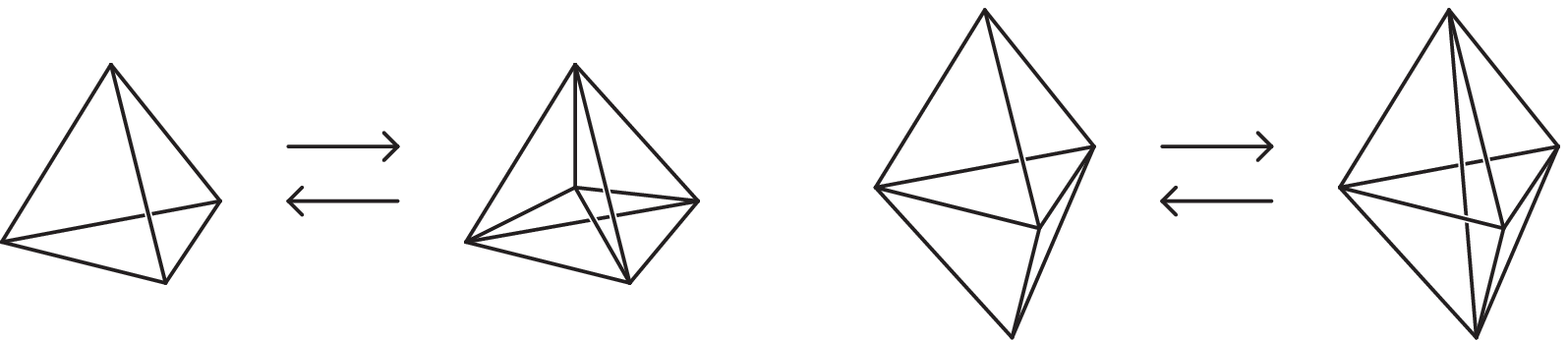}
}
\vskip 6pt
\centerline{Figure 4: Pachner moves}

Therefore, just as argued in the previous subsection, if two triangulated manifolds
are PL-homeomorphic, then one will always be able to prove this eventually.
Moreover, if one has a computable upper bound on the number of Pachner
moves required to pass between two triangulations, then one can solve
the PL-homeomorphism problem. In the case of knot exteriors, such a computable
upper bound was found by Mijatovic [71]. Unfortunately, his bound is massive: a tower of exponentials,
exponentially high. So, again, this does not lead to a practical algorithm.

Mijatovic's method built on work of Haken. As explained in Section 2.1,
Haken used a hierarchy for a link exterior $M$ to build a `canonical' cell structure for $M$.
This had the property that one can algorithmically build this cellulation starting from any
triangulation of $M$. By carefully modifying the initial triangulation to one derived
from this cellulation, Mijatovic was able to produce an upper bound on the number
of Pachner moves required. The method is not easy, and in particular, he
needed to go beyond what Haken achieved, because he also needed
the Rubinstein-Thompson machinery for 3-sphere recognition [99, 70].

\vskip 6pt
\noindent {\caps 2.6. Solving the equivalence problem efficiently}
\vskip 6pt

None of these approaches to the equivalence problem is known to be efficient.
Indeed, any efficient solution seems to be out of reach at present. This leads
to our first unsolved problem.

\noindent {\caps Problem 1.} {\sl What is the complexity class of the equivalence problem for knots and links?}

It is quite striking that, of the five above approaches to the equivalence problem,
only the final two give an {\sl a priori} estimate for its running time. 
However, the disadvantage of the approaches using Reidemeister moves or Pachner moves is
that they are almost inevitably lengthy.

One the main difficulties in Haken's approach to the equivalence problem is
that, in order to apply normal surface theory, one needs to build a triangulation
$T_i$ for each 3-manifold $M_i$ in the hierarchy. However, the surface $S_i$
may be exponentially complicated in $T_i$, in the sense of having exponentially many
triangles and squares, as a function of the number of tetrahedra in $T_i$. So, it seems inevitable that the
next triangulation $T_{i+1}$ should be exponentially more complicated than
that of $T_i$. So, the number
of tetrahedra in $T_n$ is massive. As a function of the crossing number
of a given diagram of the link, it is a tower of exponentials. The height of this
tower is $n$, the length of the hierarchy, and unfortunately, due to the
technical requirements of Haken's hierarchies, $n$ can be quite large too.

On the other hand, the approach to the equivalence using geometric structures
is of unknown complexity. However, it appears to be computationally efficient in practice,
even if one cannot yet prove this. So, there remains some hope that a provably efficient
solution to the equivalence problem may be found in the future.

One might ask for a polynomial-time algorithm in Problem 1. But this seems rather unlikely, and there certainly seems
no prospect of a proof in sight. Even an NP algorithm seems beyond reach at the moment.
(For a very brief explanation of NP and other complexity classes, see Section 3.1.)
Alternatively, one might try to find lower bounds on the complexity of the
problem, conditional upon well-known conjectures in theoretical computer science. For example,
is the equivalence problem NP-hard? In other words, is it at least as hard as any other
NP problem? In particular, if there was a polynomial-time solution to
the equivalence problem for knots and links, would this imply that P $=$ NP? This seems plausible.
Indeed, a seemingly simpler algorithm in 3-manifold theory is known to be NP-complete,
by work of Agol, Hass and Thurston [3].
This is the problem of determining whether a given knot in a given 3-manifold bounds a compact
orientable embedded surface with genus at most a given integer.

\vskip 18pt
\centerline {\caps 3. The recognition problem}
\vskip 6pt

An algorithmic problem that is closely related to the equivalence problem is the {\sl recognition problem}
for a fixed link type. Here, one fixes a link type $K$, and one asks whether a given diagram
represents this link type. Fairly obviously, a solution to the equivalence problem implies a solution
to the recognition problem for each link type. However, by fixing the link type, we simplify the
problem, and so there may be some hope of coming up with a more efficient solution.
The case of the unknot is particularly intriguing. The following remains unknown.

\noindent {\caps Problem 2.} {\sl Is there a polynomial-time algorithm to recognise the unknot?}

Currently, the unknot recognition problem is known to lie in the following complexity classes:
NP, co-NP and E. We briefly remind the reader of the definitions of these and various other
classes.

\vfill\eject
\noindent {\caps 3.1. Some basic complexity classes}
\vskip 6pt

The class P is possibly simplest to describe. It consists of those problems that be
solved in time $n^k$, where $n$ is the size of the input and $k$ is a fixed constant.
Similarly, a problem is in the class E if it can be
solved in time $k^n$, where $n$ is the size of the input and $k$ is a fixed constant.
Somewhat confusingly, this is not the same as the class EXP of exponential-time algorithms,
since, by definition, these can solved in time $k^{n^c}$ for constants $k$ and $c$.

A problem is in NP if it admits a polynomial-time certificate. This is an extra piece of
information that is not provided by the algorithm, rather it is given by some external source. The point is
that {\sl if} the answer to the problem is `yes', then the certificate proves that it is `yes',
and this can be verified in polynomial time. There are many examples of NP problems;
one is the problem of deciding whether a given positive integer is composite.
In this case, a very simple certificate is two integers (greater than 1) whose product is the given
integer. It can be verified in polynomial time (as a function of the number of digits of the input)
that these two integers do indeed multiply together to produce the given number.
One might legitimately wonder why NP problems are solvable at all, given that they
require information from an external source. But the point is that this information must have
polynomially bounded size, because a computer has time enough only to check this
many bits of information. So, a deterministic algorithm proceeds by checking all possible
certificates of at most this size. Thus, problems in NP are also in EXP.

A computational problem is in co-NP if its negation is in NP. Thus, for example, the problem of deciding
whether a positive integer is composite is in co-NP, because there is an NP algorithm to determine whether a positive integer is prime.
This is not at all obvious however. It was first proved by Pratt [84], using the Lucas primality test.

There is also the well-known class of NP-complete problems. A problem is NP-complete if it is in NP
and any other NP problem may be (efficiently) reduced to it. In particular, if an NP-complete
problem has a polynomial-time solution, then this would imply that P = NP. It is widely believed that P $\not=$ NP, 
and so NP-complete problems are, in some conditional sense, provably hard. However, it is very
unlikely that any problem that is both in NP and co-NP is NP-complete, because this would
imply that NP = co-NP, which also is viewed as rather unlikely.

\vskip 6pt
\noindent {\caps 3.2. Haken's unknot recognition algorithm}
\vskip 6pt

The first person to provide an unknot recognition algorithm was Haken [24]. A knot $K$ is the unknot
if and only if its exterior $M$ has compressible boundary. Haken showed how one may use
normal surface theory to determine this. As in Section 2.1, one builds a triangulation of 
$M$ using the given diagram of $K$. Haken showed that if $M$ contains an essential properly embedded disc,
then it contains one that is a fundamental normal surface. He also showed that all fundamental
normal surfaces may be algorithmically constructed. Thus, his unknot recognition algorithm proceeds
by constructing all possible fundamental normal surfaces, and determining whether any of these
is an essential disc.

\vfill\eject
\noindent {\caps 3.3. The NP algorithm of Hass, Lagarias and Pippenger}
\vskip 6pt

Hass, Lagarias and Pippenger [27] provided an NP algorithm, by building on Haken's work. Instead of using fundamental
surfaces, they used vertex surfaces. By definition, a normal surface $S$ is a
{\sl vertex surface} if it is connected and whenever $n [S] = [S_1] + [S_2]$ for a positive integer $n$
and normal surfaces $S_1$ and $S_2$, then each of $[S_1]$ and $[S_2]$ is a multiple
of $[S]$. The reason for the vertex terminology is that the set of all normal surface vectors
can be viewed as integer points within a larger subset of ${\Bbb R}^{7t}$,  where $t$ is the number
of tetrahedra of the 3-manifold $M$. This subset is 
the set of vectors satisfying the matching equations and quadrilateral constraints,
but where the vectors are required only to be non-negative {\sl real} numbers.
It is not hard to see that this subset is a union of convex polytopes, each of which is
the cone on a compact polytope, where the cone point is the origin.
Vertex surfaces are precisely those connected normal surfaces that
are a multiple of a vertex of one of these compact polytopes. Hass, Lagarias and Pippenger
showed that if $M$ contains an essential properly embedded disc, then it has one
that is a vertex normal surface. Moreover, for each such vertex surface that is an essential
disc, one can certify that it is an essential disc in polynomial time. Thus, unknot recognition lies in NP.
It is also not hard to prove that the number of vertex surfaces is at most $2^{7t}$,
because each one is obtained as the unique solution (up to scaling) of the matching
equations plus some extra constraints that force certain co-ordinates to be zero.
So, if one checks each such vertex surface in turn, one can determine whether the
given knot is the unknot in at most $k^n$ steps where $n$ is the crossing number of the diagram
and $k$ is a fixed constant. This implies that unknot recognition lies in E.

\vskip 6pt
\noindent {\caps 3.4. The co-NP algorithm of Kuperberg}
\vskip 6pt

Kuperberg has recently proved that unknot recognition is in co-NP, assuming the
Generalised Riemann Hypothesis [48]. Here, one wants an efficient way of certifying
that a non-trivial knot $K$ is indeed knotted. Kuperberg achieves this by establishing
the existence of a representation from $\pi_1(S^3 - K)$ to some finite group 
with non-abelian image. The unknot clearly admits no such representation,
because the fundamental group of its complement is ${\Bbb Z}$. Thus, this
does establish that the knot is non-trivial.

The existence of such a representation is not at all clear. Moreover, Kuperberg requires
a representation with small enough image, so that the representation can be
verified as a homomorphism with non-abelian image in polynomial time.
The starting point is the major theorem of Kronheimer and Mrowka [45]
that establishes that, for a non-trivial knot $K$, there is a homomorphism
$\pi_1(S^3 - K) \rightarrow SU(2)$ with non-abelian image. This is established
using a highly sophisticated argument that uses taut foliations, contact structures,
symplectic fillings, and the instanton equations. (In the case where the knot $K$ is
hyperbolic, one can instead use the associated representation $\pi_1(S^3 - K) \rightarrow PSL(2, {\Bbb C})$
in Kuperberg's argument.)

Once one has a linear representation of $\Gamma = \pi_1(S^3 - K)$, there is a standard method of obtaining
finite representations, that goes back to the result of Malce'ev [62] which states that
finitely generated linear groups are residually finite. One considers a finite generating
set for $\Gamma$, and for each generator, one considers the entries of the corresponding
matrix (in $SU(2)$). The ring $R$ generated by these entries is a finitely generated subring 
of ${\Bbb C}$. Thus, the image of $\Gamma$ lies in $SL(2,R)$. Such a ring $R$ contains many finite index ideals $I$.
The desired finite quotients are obtained by mapping $SL(2,R)$ onto $SL(2,R/I)$ for some such $I$.
The Generalised Riemann Hypothesis is used to prove the existence of an ideal $I$ with small index,
and hence one obtains a finite quotient of $\Gamma$ with small size. The actual tools that
are used here are results of Korain [43], Lagarias-Odlyzko [59] and Weinberger [104] which
imply that, for any integer polynomial, one may reduce its coefficients modulo some 
small prime and obtain a polynomial with a root, assuming the Generalised Riemann Hypothesis.

\vskip 6pt
\noindent {\caps 3.5. The co-NP algorithm of Agol}
\vskip 6pt

Agol has announced the existence of another certificate for establishing that a non-trivial knot is knotted [1].
This has the advantage that it is unconditional on any conjectures, and it also
provides more information about the knot, by determining its genus.
Recall that the {\sl genus} of a knot $K$ is the minimal genus of a Seifert
surface for $K$. It is a fundamental quantity, but for the purposes of this algorithm,
all that one needs is that the genus is zero if and only if the knot is trivial.

It is a consequence of Haken's work that the genus of a knot is algorithmically computable,
because a minimal genus Seifert surface can be arranged to be a fundamental normal surface.
However, Thurston and Gabai found another method for determining the genus of knots, by using the
theory of taut foliations [101, 21]. In fact, this method naturally measures not the genus of a compact surface $S$
but its {\sl Thurston complexity} $\chi_-(S)$. By definition, when $S$ is connected,
then $\chi_-(S) = \max \{ -\chi(S), 0 \}$. When $S$ is disconnected with
components $S_1, \dots, S_n$, then $\chi_-(S) = \sum_{i=1}^n \chi_-(S_i)$.
Thus, when $S$ is orientable, $\chi_-(S)$ is roughly the negative of the
Euler characteristic of the surface, but one first discards disc and sphere components.
So, for knots, a Seifert surface has minimal possible $\chi_-$ if and only if it has
minimal genus. But for links, this need not be the case, because the link may have
disconnected Seifert surfaces as well as connected ones.

The theory of Thurston and Gabai is concerned with foliations on a 3-manifold, where the leaves of the foliation
have codimension one and where there is a consistent transverse orientation on these leaves.
Such a foliation is {\sl taut} if there is a simple closed curve in
the 3-manifold that is transverse to the foliation and that intersects every leaf. This may
seem a slightly strange definition, but it turns out that the existence of such a curve
rules out various trivial examples of foliations, and in fact the following was proved by Thurston [101]:
{\sl If $S$ is a compact leaf of a taut foliation on a compact orientable 3-manifold $M$,
then $S$ has minimal $\chi_-$ in its homology class in $H_2(M, \partial S)$.}
Hence, if one starts with a Seifert surface $S$ for an oriented link $L$, and one can find a taut
foliation on the exterior of $L$ in which $S$ is a compact leaf, then $S$ has minimal
$\chi_-$ among all Seifert surfaces for $L$. Gabai [21] introduced a method for constructing
these foliations and in fact proved the converse statement: {\sl If $S$ is a Seifert
surface for an oriented non-split link $L$ with minimal $\chi_-$, then there is some
taut foliation on the exterior of $L$ in which $S$ is a leaf.}

One can view this taut foliation as a way of certifying the genus of a knot. However,
it is not a certificate in the strict algorithmic sense, because in principle, an infinite
amount of information is required to specify the foliation. But the key features
of Gabai's foliations are encoded by a finite amount of information as follows.
Gabai utilised hierarchies to construct his foliations, and it was first observed
by Scharlemann [88] that many of the useful topological consequences of Gabai's
theory can be extracted solely from the hierarchy, without referring to the foliation.
The hierarchies of Gabai and Scharlemann naturally occur in the setting of 3-manifolds
with the following extra structure.

A {\sl sutured manifold} is a compact oriented 3-manifold $M$ with its boundary
decomposed into two subsurfaces $R_-$ and $R_+$. A transverse orientation
is imposed on these subsurfaces, so that $R_+$ points out of $M$ and
$R_-$ points into it. These surfaces intersect along a collection of simple closed curves
$\gamma$, called {\sl sutures}. It is usual to denote a sutured manifold by $(M, \gamma)$,
but the surfaces $R_+$ and $R_-$ are also an essential part of the structure. A sutured manifold 
$(M, \gamma)$ is {\sl taut} if $M$ is irreducible, 
$R_+$ and $R_-$ are incompressible and they minimise $\chi_-$ in their homology classes in $H_2(M, \gamma)$.

When a sutured manifold $(M,\gamma)$ is decomposed along a properly embedded,
transversely oriented surface $S$ that intersects $\gamma$ transversely, the resulting 3-manifold $M'$ inherits a sutured manifold
structure. This decomposition is {\sl taut} if $M$ and $M'$ are both taut as sutured manifolds
and $S$ has minimal $\chi_-$ in its homology class in $H_2(M, \partial S)$.
A sequence of such decompositions terminating in a collection of 3-balls is
a {\sl sutured manifold hierarchy}. A key feature of sutured manifold theory is:
{\sl If a sequence of sutured manifold decompositions ends in a collection of
taut sutured 3-balls, then (provided some simple conditions hold)
every decomposition in this sequence is taut}. 

This sequence of decompositions can therefore be viewed as a certificate for knot genus.
However, it is not at all clear that one can produce such a certificate that can be verified in 
polynomial time. The principal difficulty is the one explained in Section 2.6 in reference
to Haken's use of hierarchies. It seems, at first sight, that one needs to keep track of
a triangulation of each manifold in this hierarchy, and the complexity of these triangulations
seems to grow rapidly as one proceeds along the hierarchy. Agol's alternative approach
is to make the entire sutured manifold hierarchy normal (in some suitable sense) with the respect
to the initial triangulation of $M$. The idea is that the sutured manifold hierarchy is closely
related to a taut foliation and so normalising such a hierarchy is somewhat analogous to
normalising a taut foliation. Techniques exist for placing taut foliations in `normal' form, due to
Brittenham [8] and Gabai [22], and versions of these are used by Agol. Unfortunately, full
details of this proof have not yet been written down. An alternative approach may be
to use the techniques of the author in [50], where sutured manifold hierarchies were normalised
in some weak sense.\footnote{$^2$}{\ninerm This approach has been successfully completed
by the author [56], giving the first full unconditional proof that unknot recognition is in co-NP.}

\vskip 6pt
\noindent {\caps 3.6. Khovanov homology}
\vskip 6pt

In the mid 1980s, knot theory underwent a dramatic revolution, with the
introduction of the Jones polynomial [37], and then Witten's interpretation of this using
Chern-Simons theory [105]. Right from the outset of this work, it was asked:
{\sl Does the Jones polynomial detect the unknot?} In other words, must a
non-trivial knot necessarily have non-trivial Jones polynomial? This question remains
unanswered. Indeed, it may be viewed as a central open problem in knot theory,
but it is not `elementary' and so does not make it onto our problem list.

It was shown recently by Kronheimer and Mrowka [47] that a more refined invariant, Khovanov homology,
does detect the unknot. This followed from a slightly indirect relationship
between Khovanov homology and 4-dimensional gauge-theoretic invariants.

It seems unlikely that this result will lead to an efficient way of
recognising the unknot. This is because determining Khovanov homology
is a computationally intensive process. In fact, the less refined Jones polynomial
is, in one sense, provably hard to calculate. For example, it was shown by 
Jaeger, Vertigana and Welsh [35] that
computation of the Jones polynomial of a link is `$\sharp$ P-hard'.
(This is somewhat similar to being NP-hard.)

\vskip 6pt
\noindent {\caps 3.7. Heegaard Floer homology}
\vskip 6pt

Another major advance in knot theory has been the introduction of Heegaard Floer homology
by Oszv\'ath and Szab\'o [77].
This was based on earlier instanton invariants, as pioneered by Donaldson [16]
in dimension 4, and then Floer [19] in dimension 3, and monopole
invariants, as developed by Seiberg and Witten [106] in dimension 4 and then Kronheimer and
Mrowka [46] in dimension 3. It was clear from the outset that
these gauge-theoretic invariants are particularly powerful, but they were undoubtedly hard
to calculate. One of the most attractive features of Heegaard Floer homology is that,
right from its inception, it has been possible to calculate it in many important examples. 
Indeed, there are now algorithms to calculate the various different versions of
Heegaard Floer homology, due to Sarkar and Wang [86]
and Manolescu, Ozsv\'ath and Sarkar [64].
One of the earliest results about Heegaard Floer homology was that it detects
the unknot [78]. So, coupled with the fact that it is computable,
this provides another unknot recognition algorithm. Again, however, it
seems unlikely that this algorithm is in any way efficient.

\vskip 6pt
\noindent {\caps 3.8. Rectangular diagrams and arc presentations}
\vskip 6pt

In 2003, Dynnikov introduced a new and striking solution to the unknot recognition problem [17].
This used rectangular diagrams and arc presentations. A {\sl rectangular diagram}
(also called {\sl grid diagram}) is a planar link diagram that consists of vertical
and horizontal arcs, no two of which are colinear, and with the key property that
whenever a horizontal and a vertical arc cross, it is the vertical arc that is the over-arc
at the resulting crossing. The number of vertical arcs is equal to the number of
horizontal arcs; this is known as the {\sl arc index} of the rectangular diagram.

\vskip 18pt
\centerline{
\includegraphics[width=1.5in]{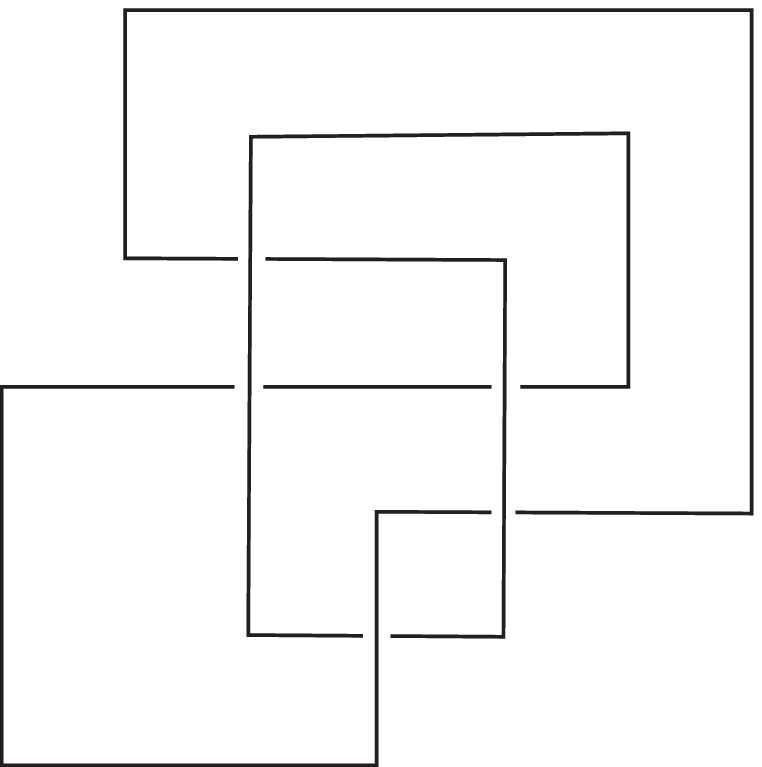}
}
\vskip 6pt
\centerline{Figure 5: A rectangular diagram}

Rectangular diagrams were studied in detail by Cromwell [12]. He introduced a set of
moves which change a rectangular diagram, without changing the link type, much like
Reidemeister moves. He showed that any two rectangular diagrams for a link differ
by a sequence of such moves. (See Figure 6.) Two of these types of moves ({\sl cyclic permutations}
and {\sl exchange moves}) leave the arc index unchanged. The other type of move
is a {\sl stabilisation/destabilisation} which increases/decreases the arc index by one.

\vskip 18pt
\centerline{
\includegraphics[width=5.4in]{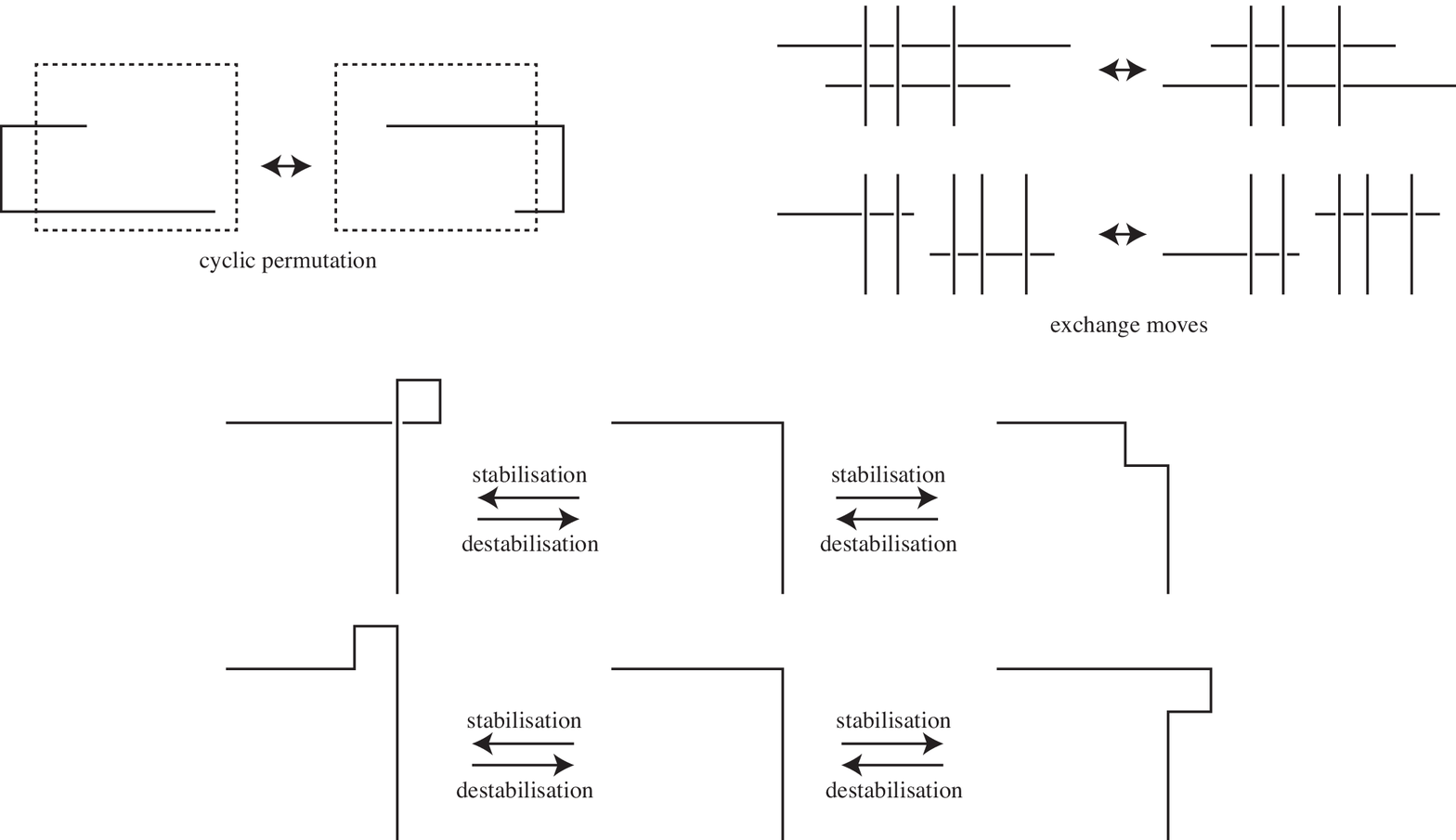}
}
\vskip 6pt
\centerline{Figure 6: Modifications to a rectangular diagram}

Dynnikov proved that any rectangular diagram of the unknot can be reduced to the
trivial diagram using these moves, but crucially, {\sl no stabilisations are required}.
This leads to the following unknot recognition algorithm. Start with the given diagram,
and using an isotopy of the plane, make it rectangular. Then keep applying all possible cyclic
permutations, exchange moves and destabilisations. Only finitely many moves
are possible at each stage. The key point is that the number of rectangular
diagrams with bounded arc index is finite. So, eventually, one can explore the
entire set of rectangular diagrams that are reachable from the initial one by 
sequences of these moves. The given knot is the unknot if and only if the trivial
diagram is in this collection.

This does not immediately lead to an efficient solution. This is because the
number of rectangular diagrams with a given arc index is a super-exponential
function, and so the search space might be quite large. But it is a beautiful 
and striking result that is perhaps the closest in spirit to elementary knot theory
among all the known unknot recognition algorithms. It also has the following
consequence, as noted by Dynnikov himself: {\sl If a diagram of the unknot
has $n$ crossings, then there is a sequence of Reidemeister moves taking
it to the trivial diagram, where each diagram in this sequence has crossing number
at most $2(n+1)^2$.} This is simply because a rectangular diagram with
arc index $m$, say, has crossing number at most $(m-1)^2/2$.

Dynnikov's proof builds on previous work of Cromwell [12], Birman, Menasco [7] and Bennequin [4].
Cromwell observed a correspondence between rectangular diagrams and arc presentations,
which are defined as follows. One fixes a round unknot in the 3-sphere, which is a reference
circle, known as the {\sl binding circle}. The complement of the binding circle is an
open solid torus, which is foliated by open discs called {\sl pages}. A link is in an
{\sl arc presentation} if it intersects the binding circle in finitely many points called {\sl vertices},
and it intersects each page in either the empty set or a single open arc joining
distinct vertices. Cromwell proved that there is a one-one correspondence between
arc presentations and rectangular diagrams up to cyclic permutation. 

Thus, Dynnikov's proof mainly concentrates on arc presentations. Given an
arc presentation of the unknot, its spanning disc can be arranged so that
its intersection with the pages induces a singular foliation on the spanning disc.
The singularities occur when the disc slices through the binding circle,
and also at Morse-type singularities in the complement of the binding
circle. Dynnikov's measure of complexity for the disc is the number of
singularities of this foliation. Using an argument that relies on the fact that the disc has positive
Euler characteristic, Dynnikov showed that this singular foliation must have,
at some point, a certain configuration. He then showed that this implies
the existence of some exchange moves, cyclic permutations and possibly a destabilisation,
after which the arc index has been reduced or the number of singularities of the disc has
been reduced. Thus, this process is guaranteed to terminate. In particular, it is possible
to simplify the disc without using stabilisations. Once the disc is as simple
as possible, so too is the rectangular diagram.

This style of argument was very much based on earlier work of Cromwell,
Birman, Menasco and Bennequin. For example, Birman and Menasco studied
representations of the unknot as a closed braid, and proved that it may be
reduced to the trivial braid (of braid index one) by a sequence of moves,
also known as exchange moves and destabilisations. Thus, the braid
index does not increase during this sequence of moves. However, this does
not immediately lead to an unknot recognition algorithm, because
in the case of braids, there are infinitely many possible exchange moves
that can be applied at each stage, and the number of braids with a given
braid index is infinite. Nevertheless, Birman and Hirsch [6] did use braids
to produce an unknot recognition algorithm, by combining these singular foliation
techniques with Hass and Lagarias's work on normal surfaces [26].

\vskip 18pt
\centerline{\caps 4. Reidemeister moves}
\vskip 6pt

Reidemeister proved that any two diagrams of a link differ by a sequence of Reidemeister moves [85].
This result has many applications. For example, it is often used to show that an invariant of link diagrams
actually leads to an invariant of knots and links, by showing that the invariant is unchanged
under each Reidemeister move. It also leads to many natural and interesting questions, including the
following.

\vskip 6pt
\noindent {\caps Problem 3.} {\sl Find good upper and lower bounds on the number of Reidemeister
moves required to relate two diagrams of a knot or link.}
\vskip 6pt

In addition to being natural and attractive, this question has algorithmic implications. As explained in
Section 2.4, there exists a computable upper bound on Reidemeister moves if and only if
there exists a solution to the equivalence problem. Moreover, if one just considers diagrams
of a fixed link type, then the existence of a computable upper bound on Reidemeister moves
is equivalent to the existence of a solution to the recognition problem for that link type.
Obviously, one would like upper bounds that are as small as possible, because this
leads to more efficient algorithms.

We now give a summary of what is known here. The gap between the known upper and lower bounds are quite startling.

\vskip 6pt
\noindent {\caps 4.1. Lower bounds}
\vskip 6pt

There is, of course, a linear lower bound on the number of Reidemeister moves required to relate two diagrams of
a link. More precisely, if two diagrams have $n_1$ and $n_2$ crossings, then the number of moves
relating them is at least $|n_1 - n_2|/2$, simply because each Reidemeister move changes the crossing
number by at most 2.

Hass and Nowik [28] proved that, in general, at least quadratically many moves may be needed. For each link type,
they produced a sequence of diagrams $D_n$, where the crossing number of $D_n$ grows
linearly with $n$, but where the minimal number of Reidemeister moves relating $D_n$ to $D_0$ is
a quadratic function of $n$. They did this by finding an invariant of knot diagrams that does
not lead to a knot invariant, but that changes in a controlled way when
a Reidemeister move is performed. Hence, a `large' difference in the invariants of two diagrams
implies that they necessarily differ by a long sequence of Reidemeister moves.

\vskip 6pt
\noindent {\caps 4.2. Upper bounds for general knots and links}
\vskip 6pt

In contrast to these lower bounds, the known upper bounds on Reidemeister moves
are vast. The best known bound in general is due to Coward and the author [9]:
{\sl If two diagrams of a link have $n$ and $n'$ crossings
respectively, then these diagrams differ by a sequence of at most 
$$2^{2^{\cdot^{\cdot^{\cdot^{2^{(n+n')}}}}}} \Big\}\hbox{height } c^{(n+n')}$$
Reidemeister moves, where $c = 10^{1000000}$.}
This is huge, but it was the first known upper bound that is primitive recursive.

The proof relied on Mijatovic's upper bound on Pachner moves between triangulations
of a link exterior, as described in Section 2.5,
which in turn built on Haken's work.
Mijatovic's upper bound was a tower of exponentials, much like the 
formula above. However, it is not at all clear how one can go from a bound on Pachner
moves to a bound on Reidemeister moves. The idea is to use the sequence of Pachner moves
to build an explicit homeomorphism from the 3-sphere to itself taking the link
arising from the first diagram to the link arising from the second diagram. This homeomorphism
is isotopic to the identity (via the `Alexander trick'), and so one obtains a one-parameter family of links interpolating
between the two given ones. Projecting these links gives a one-parameter family of
diagrams. The control that we have over the homeomorphism of the 3-sphere
can be used to control the complexity of these link projections, and so one obtains an upper bound
on the crossing number of all the intermediate diagrams. This in turn can be used to
bound the number of Reidemeister moves (as explained in more detail in Section 4.5).

\vfill\eject
\noindent {\caps 4.3. Upper bounds for the unknot}
\vskip 6pt

A recent theorem of the author [55] provides a polynomial upper bound on Reidemeister moves for the unknot:
{\sl If a diagram of the unknot has $n$ crossings, then there is a sequence of at
most $(236n)^{11}$ Reidemeister moves taking it to the trivial diagram.}

Previously, the best known upper bound was an exponential function of $n$. This was due
to Hass and Lagarias [26]. Their method of proof relied heavily on Haken's solution of the recognition problem for the unknot.
One is given a diagram of the unknot with $n$ crossings, and from this one builds
a triangulation of the unknot's exterior. It is not hard to arrange that each tetrahedron
of this triangulation is straight in ${\Bbb R}^3$ (apart from the tetrahedron enclosing the 
point at infinity) and that the number of tetrahedra is bounded
above by a linear function of $n$. Now a spanning disc
for the unknot can be placed in normal form, and can in fact be realised as a vertex
surface. Crucially, there is an upper bound on the number of triangles and
squares in a vertex normal surface, which is an exponential function of the number of tetrahedra
in the triangulation. One then simply slides the unknot along this disc.
Each time one slides across a triangle or square, this induces a controlled
number of Reidemeister moves in the projection to the plane of the diagram.

So, it is the exponential upper bound on the size of the normal spanning disc
that leads to the Hass-Lagarias upper bound on the number of Reidemeister moves.
It is possible to prove that the exponential upper bound on the size of the spanning
disc cannot be improved upon. Hass, Snoeyink and Thurston [29] found a sequence
of unknot diagrams $D_n$, where the crossing number is a linear function of $n$,
but where any piecewise-linear spanning disc for the unknot must consist of
at least exponentially many triangles.

To prove that there is a polynomial upper bound on Reidemeister moves, the
author had to use both normal surfaces and arc presentations. As in Dynnikov's
work (which is described in Section 3.8), one starts with an arc presentation of
the unknot, and the aim is to apply exchange moves, cyclic permutations and
destabilisations to reduce it to the trivial arc presentation. Dynnikov examined 
a spanning disc and its induced singular foliation, and showed that there is
always a move that can be performed which results in a reduction in the
number of singularities in the foliation or a reduction in the arc index. Using normal surface theory, it is
not hard to show that initially, the disc can be arranged to have at most exponentially
many singularities. So, although this gives an upper bound on the number of
moves, this bound is exponential. To obtain an improved estimate, one
needs to go deeper into normal surface theory. It is possible to speak of
two bits of normal surface as being {\sl normally parallel}. The aim is to find
many parallel pieces of the disc which can be moved all at the same time.
In this way, one obtains, using a controlled number of exchange moves
cyclic permutations and destabilisations, a substantial decrease in the complexity of the disc.
This decrease is large enough that only polynomially many moves are
required before the disc has been completely simplified, which implies that
the arc presentation has become trivial.

\vskip 6pt
\noindent {\caps 4.4. An upper bound for each link type}
\vskip 6pt

The polynomial upper bound on Reidemeister moves described above has recently
been generalised by the author, as follows [58]. {\sl For each link type $K$, there is a polynomial
$p_K$ such that any two diagrams of $K$ with crossing number $n$ and $n'$
differ by a sequence of at most $p_K(n) + p_K(n')$ Reidemeister moves}.
As  a consequence, {\sl for each link type $K$, the $K$-recognition problem is in NP.}
The certificate is simply a sequence of Reidemeister moves taking the given diagram
of $K$ to some fixed diagram.

The proof uses and extends the techniques that were developed in the case of the unknot.
In particular, both normal surfaces and arc presentations play a central role. In the case
of the unknot, the spanning disc was simplified at each stage of the procedure.
For a general knot $K$, one uses instead an entire hierarchy for the exterior of $K$.
The spanning disc for the unknot could be simplified until its singular foliation
contained just two singularities. This is not possible for the surfaces in a hierarchy.
Instead, the goal is to simplify these surfaces until their number of singularities is
bounded by a polynomial function of the initial arc index.  Once this has been achieved, one considers
a regular neighbourhood of this hierarchy together with a regular neighbourhood of the knot.
This is a ball $B$, and the complexity of $B$ is, in some suitable sense,
controlled. The ball inherits a handle structure by thickening the cell structure of the hierarchy, and $K$
lies inside this handle structure in a way that depends only on the topology of the hierarchy. The proof proceeds
by isotoping $K$ through $B$ until it lies within a single 0-handle and projects
to some fixed diagram for $K$. Since $B$ has controlled complexity, so too
does the projection of $K$ throughout this isotopy. Thus, one obtains a bound
on the crossing number of the intermediate diagrams, and with a bit more work,
a polynomial bound on the number of Reidemeister moves used in this process.

\vskip 6pt
\noindent {\caps 4.5. Bounds on the crossing number of intermediate diagrams}
\vskip 6pt

In addition to bounding the number of Reidemeister moves required to pass between two link
diagrams, one can also seek to bound the complexity of the diagrams in this sequence.

\vskip 6pt
\noindent {\caps Problem 4.} {\sl Given two diagrams of a link, find a sequence of Reidemeister moves
relating them where all diagrams in the sequence have small crossing number.}
\vskip 6pt

This is, of course, related to Problem 3. Firstly, an upper bound on the number of moves
immediately gives an upper bound on the intermediate crossing numbers. This is because
each Reidemeister move changes the crossing number by at most 2. Secondly, if one has an upper
bound on the crossing number of each intermediate diagram, then one obtains an upper bound on
the number of moves. This is because there are at most $(24)^{n+1}$ connected diagrams with crossing number
at most $n$ (see Theorem 6.5 in [9] for example). Moreover, in a shortest sequence of Reidemeister moves joining two diagrams,
one never visits the same diagram twice. So, if each diagram in the sequence has crossing number
at most $n$, say, then the length of the sequence is at most $(24)^{n+1}$.

However, one can sometimes establish a much better bound than this. As explained in Section 3.8,
Dynnikov found, in the case of the unknot, a quadratic upper bound on the crossing number
of intermediate diagrams.

\vfill\eject
\centerline{\caps 5. Crossing number}
\vskip 6pt

The {\sl crossing number} $c(K)$ of a knot or link $K$ is defined to be the minimal number
of crossings in any diagram for $K$. 

\vskip 6pt
\noindent {\caps 5.1. Computing the crossing number}
\vskip 6pt

Crossing number is algorithmically computable. Indeed, this fact is closely related
to the solubility of the equivalence problem for links. For example, the problem of determining
whether a knot's crossing number is zero is clearly equivalent to determining whether it is
the unknot. One can use a solution to the equivalence problem to determine a knot's
crossing number in the following naive way. Starting with a knot diagram having $n$ crossings,
enumerate all diagrams with less than $n$ crossings, and for each, determine whether it
is the same knot as the original one. The one with the smallest number of crossings clearly
realises the knot's crossing number. This approach is obviously not very efficient but
it seems unlikely that there is any quicker way of determining a link's crossing number in general.
In fact, this method has been successfully used in practice to 
compile knot tables, by Hoste, Thistlethwaite and Weeks for example [32].

\vskip 6pt
\noindent {\caps 5.2. Determining the crossing number for certain link classes}
\vskip 6pt

Although there seems little hope of anything more than using a brute-force search to determine
the crossing number of an arbitrary link, one can often determine the crossing number more efficiently in various interesting cases. For example, the crossing numbers of alternating links are completely understood,
by work of Kauffman [40], Murasugi [72] and Thistlethwaite [97]. Their work, and a broader analysis
of alternating links, is discussed in Section 7.1.

Torus links form another interesting collection, although they are of course more restricted.
Murasugi [73] proved that the crossing number of the $(p,q)$ torus link is $pq - \max \{ p,q \}$.

\vskip 6pt
\noindent {\caps 5.3. The number of knot types}
\vskip 6pt

Although crossing number is computable for any given knot, its behaviour is
not well understood for infinite collections of knots. In particular, the following problem
remains unresolved.

\noindent {\caps Problem 5.} {\sl If $N(c)$ is the number of knot types with crossing number $c$,
determine the asymptotic behaviour of $N(c)$ as $c$ tends to infinity.}

There is a more-or-less obvious exponential upper bound to $N(c)$. As explained in Section 4.5,
the number of connected diagrams with crossing number $c$ is at most $24^{c+1}$, and so $N(c) \leq 24^{c+1}$.

One can also obtain an exponential lower bound on $N(c)$, by considering alternating knots.
As described in the previous subsection, the crossing numbers of alternating knots are completely
understood: a `reduced' alternating diagram has minimal crossing number. (For an explanation
of this terminology, see Section 7.1.) Moreover, the resolution of the Tait Flyping conjecture
by Menasco and Thistlethwaite [69] implies that one can determine exactly when two alternating
diagrams represent the same knot type. So, enumerating the number of alternating knot
types with crossing number $c$ is nearly the same as enumerating the number of
alternating diagrams. Very accurate asymptotics were established by Sundberg and
Thistlethwaite [95] who showed that, if $N_{\rm alt}(c)$ is the number of prime alternating link types
with crossing number $c$, then
$$\lim_{c \rightarrow \infty} (N_{\rm alt}(c))^{1/c} = {101 + \sqrt{21001} \over 40}.$$

Using similar methods, Thistlethwaite [98] was also able to show that alternating links are increasingly scarce 
among all link types, as the crossing number tends to infinity. In particular, $N_{\rm alt}(c) / N(c) \rightarrow 0$
as $c \rightarrow \infty$. But more precise asymptotics for $N(c)$ seem a long way from resolution.

\vskip 6pt
\noindent {\caps 5.4. The crossing number of composite knots}
\vskip 6pt

One of the most basic operations in knot theory is connected sum. Here, one starts with two oriented knots $K_1$ and $K_2$
in the 3-sphere. For each knot, one finds a 3-ball that intersects the knot in a properly embedded unknotted arc.
One removes the interior of this ball. The result is, for each knot, a knotted arc in the complementary 3-ball.
If one glues these two balls together, the result is the 3-sphere. One performs this gluing so that the two arcs
join together to form a knot, and so that the orientations on the arcs patch together correctly. The resulting knot
is the {\sl connected sum} of $K_1$ and $K_2$, denoted $K_1 \sharp K_2$.

\vskip 18pt
\centerline{
\includegraphics[width=4in]{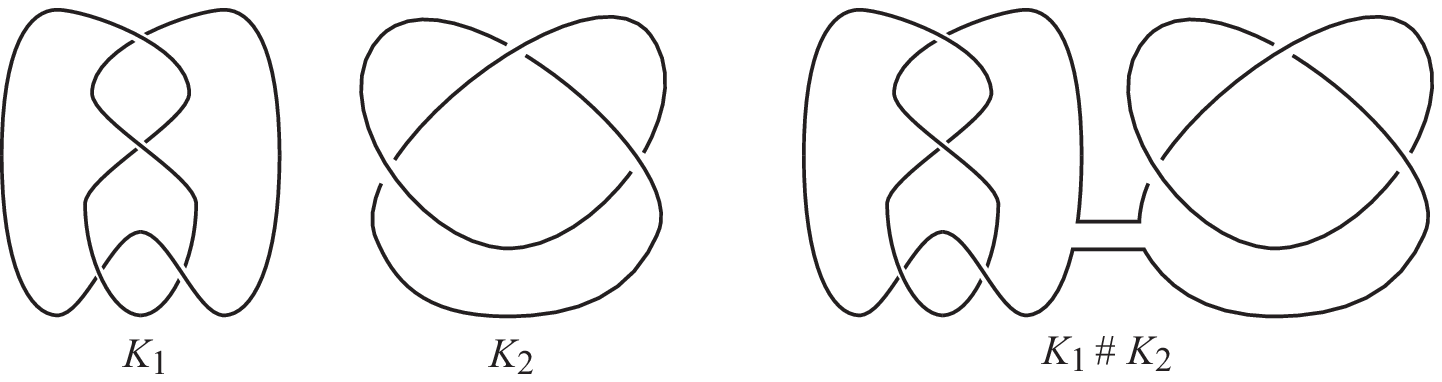}
}
\vskip 6pt
\centerline{Figure 7: Connected sum}

As is evident from Figure 7, one can create a diagram $D$ for $K_1 \sharp K_2$ starting from diagrams of
$K_1$ and $K_2$. If we use diagrams of $K_1$ and $K_2$ with minimal crossing number, we obtain
a diagram $D$ with crossing number $c(D)$ equal to $c(K_1) + c(K_2)$. Hence, $c(K_1 \sharp K_2) \leq c(K_1) + c(K_2)$.
It seems `obvious' that this should be an equality, but actually, it is far from clear, and this is in fact
one of the most notorious unsolved problems in elementary knot theory.

\noindent {\caps Problem 6.} {\sl Must $c(K_1 \sharp K_2) = c(K_1) + c(K_2)$?}

This has been verified for certain knot classes, for example when $K_1$ and $K_2$ are both
alternating [40, 72, 97] or when $K_1$ and $K_2$ are both torus knots [15]. The only general result in this
direction is due to the author [52]: {\sl For oriented knots $K_1$ and $K_2$,}
$${c(K_1) + c(K_2) \over 152} \leq c(K_1 \sharp K_2) \leq c(K_1) + c(K_2).$$

We will give an overview of the proof, because it uses some unexpected methods
which may be applicable to other problems in elementary knot theory.

There is a simple operation on knots and links
that is closely related to connected sum. This is {\sl disjoint union}.
Here one starts with two links $L_1$ and $L_2$, each lying in the 3-sphere.
One encloses each link in the interior of a ball, and then one identifies the boundaries of
these balls, thereby forming a new link in the 3-sphere, which is denoted
$L_1 \sqcup L_2$. 

Now, crossing number {\sl does} behave well with respect to disjoint union.
It is easy to prove that $c(L_1 \sqcup L_2) = c(L_1) + c(L_2)$. Again we have the upper
bound $c(L_1 \sqcup L_2) \leq c(L_1) + c(L_2)$, because one can start with diagrams for
$L_1$ and $L_2$ and use these to create a diagram for $L_1 \sqcup L_2$. To prove
the inequality in the other direction, one must start with a diagram of $L_1 \sqcup L_2$
(which we may take to have minimal crossing number) and use this to create
diagrams for $L_1$ and $L_2$. But this is straightforward: to form the diagram for
$L_1$, say, one simply throws away the parts of the link projection that come from
$L_2$. Thus, one obtains diagrams for $L_1$ and $L_2$, with the property that the
sum of the number of crossings in these diagrams is at most the crossing number of the
original diagram of $L_1 \sqcup L_2$. Thus, $c(L_1) + c(L_2) \leq c(L_1 \sqcup L_2)$.

One can immediately see from this the main difficulty in proving that $c(K_1) + c(K_2) \leq c(K_1 \sharp K_2)$.
If one starts with a diagram $D$ for $K_1 \sharp K_2$ with minimal crossing number, then there is no
obvious way of forming diagrams for $K_1$ and $K_2$ from this. Although one might view a part
of the curve $K_1 \sharp K_2$ as coming from $K_1$, say, this part is only an arc and so one needs
to add in something extra to form a closed curve representing $K_1$. This may dramatically increase
the number of crossings. This is evident when one considers the 2-sphere that is used to form
the connected sum. If one performs an isotopy on $K_1 \sharp K_2$ so that it projects to give the
diagram $D$, the 2-sphere may be extremely distorted in 3-space. It is precisely an arc on
this sphere that needs to be inserted to form a copy of $K_1$. The proof of the lower bound on $c(K_1 \sharp K_2)$ proceeds by
controlling the complexity of this 2-sphere. In fact, one cuts $K_1 \sharp K_2$ at its two points of intersection
with the sphere, and then one adds two arcs, one on either side of this sphere, to form a copy of $K_1 \sqcup K_2$.
By controlling the sphere, one can control the number of new crossings this introduces. One thereby
finds a diagram for $K_1 \sqcup K_2$ with crossing number at most $152 \, c(D)$, which proves the theorem.

The way that the 2-sphere is controlled is via normal surface theory. One starts with a diagram
$D$ for $K_1 \sharp K_2$ with minimal crossing number. From this one can build a triangulation
of the exterior $M$ of $K_1 \sharp K_2$ in a reasonably natural way. The restriction of the 2-sphere
to $M$ is a properly embedded annulus $A$. Since $A$ is essential, it can be isotoped into normal form.
Hence, $A$ intersects each tetrahedron in a collection
of triangles and squares, as shown in Figure 2. Within a tetrahedron, there are four types of
triangle and three types of square. Between adjacent triangles or squares of the same type,
there lies a product region. The union of these product regions forms an $I$-bundle ${\cal B}$ which 
is embedded in the exterior of $A$. This is known as the {\sl parallelity bundle} for the exterior of $A$.
The $\partial I$-bundle is the {\sl horizontal boundary} of ${\cal B}$, denoted
$\partial_h {\cal B}$. Note that $\partial_h {\cal B}$ lies in the two copies of $A$ in
the exterior of $A$. The key part of the argument establishes that one can find properly embedded arcs
in these two copies of $A$, each joining the two boundary components of the annulus,
and which avoid $\partial_h{\cal B}$. The significance of this is that these arcs then can be made
to intersect each tetrahedron of the triangulation in a restricted way. To see this, note that
in each tetrahedron, all but at most 6 pieces in the exterior of a normal surface lie in ${\cal B}$.
Since these arcs intersect each tetrahedron in a controlled way, one can bound from above the
crossing number of their diagrammatic projection. Thus, we obtain a diagram for $K_1 \sharp K_2$,
and with some work, one can show that it has at most $152 \,c(D)$ crossings.

\vskip 6pt
\noindent {\caps 5.5. The crossing number of satellite knots}
\vskip 6pt

A natural generalisation of connected sum is the satellite knot construction.
A knot $K$ is said to be a {\sl satellite} of a non-trivial knot $L$ if $K$ lies in a regular
neighbourhood $N(L)$ of $L$, but it does not lie within a 3-ball inside $N(L)$
and is not a core curve of $N(L)$. The knot $L$ is called the {\sl companion}.

\vskip 18pt
\centerline{
\includegraphics[width=3.5in]{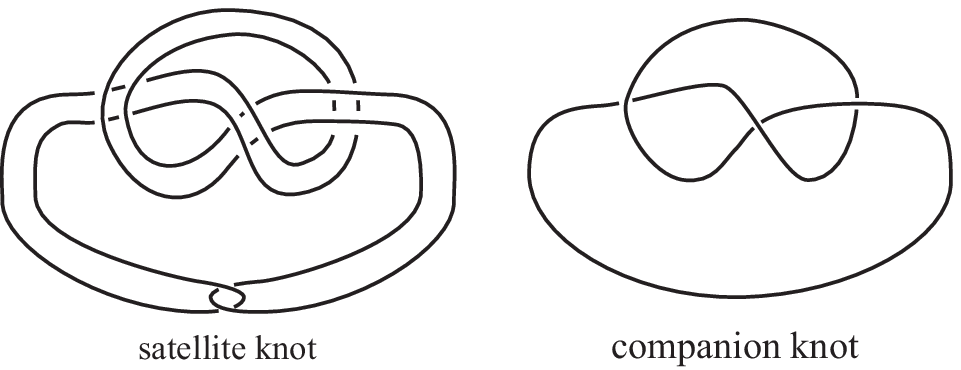}
}
\vskip 6pt
\centerline{Figure 8: A satellite knot}

The {\sl wrapping number} of the satellite is the minimal number of intersections
between $K$ and a meridian disc for $N(L)$. It is always positive.
Note that when $K_1$ and $K_2$
are non-trivial knots, then $K_1 \sharp K_2$ is a satellite of $K_1$ (and of $K_2$)
with wrapping number 1.

The following remains famously unsolved.

\noindent {\caps Problem 7.} {\sl If $K$ is a satellite of $L$, must $c(K)$ be at least $c(L)$?
More generally, must $c(K)$ be at least $w^2 c(L)$, where $w$ is the wrapping number
of the satellite?}

The motivation here is that there is a way of building a diagram for the satellite $K$ starting
from a diagram for $L$. Each crossing of $L$ becomes $w^2$ crossings for $K$. There may
be other crossings, and so unlike in Problem 6, this is a conjectured inequality rather than equality.
One might be tempted to speculate that there is an equality relating $c(K)$ to the crossing
number of $L$, the wrapping number and the crossing number of the knot in $N(L)$,
once this has been suitably defined. We do not do so here, although it is possible that
something along these lines may be true.

Unsurprisingly, the methods used to analyse the crossing number of composite knots extend to the case of satellite knots,
giving the following result of the author [54]: {\sl If $K$ is a satellite of a knot $L$, then $c(K) \geq 10^{-13} c(L)$.}

In the proof, one starts with a diagram $D$ of $K$ with minimal crossing number, and the goal
is to create a diagram of $L$ with crossing number at most $10^{13} c(D)$. One uses $D$
to construct a handle structure for the exterior of $K$. The satellite torus $T$ can be placed in `normal'
form with respect to this handle structure. Cutting along $T$ gives two pieces, one a copy of the exterior
of $L$, the other $N(L) - {\rm int}(N(K))$. Each piece inherits a handle structure.
We then modify this handle structure, primarily by replacing
parts of the parallelity bundle which are $I$-bundles over discs by 2-handles. These 2-handles may
end up lying in 3-space in a complicated way, but the remainder of the handle structure is controlled.
We then re-attach $N(K)$ to $N(L) - {\rm int}(N(K))$ to form a handle structure of the solid torus $N(L)$. A theorem of
the author [53] gives that in any handle structure of the solid torus, there is a core curve which lies
entirely within the 0-handles and 1-handles, and which intersects each such handle in a restricted way.
In this case, the core curve is a copy of $L$, which then lies inside 3-space in a controlled way. Projecting this
gives the required diagram of $L$.

\vskip18pt
\centerline{\caps 6. Crossing changes and unknotting number}
\vskip 6pt

An elementary operation that one can perform on a knot is to change a crossing in
one of its diagrams, as in Figure 9. This is strikingly difficult to analyse, and it leads to some of the most
challenging questions in the subject. An excellent survey of this topic can also be found in [90].

\vskip 6pt
\noindent {\caps 6.1. Computing unknotting number}
\vskip 6pt

The {\sl unknotting number} $u(K)$ of a knot $K$ is the minimal number of crossing changes
that one can make to some diagram that turns it into the unknot.
The difficulty of course is that one quantifies over all possible diagrams, and so although any diagram for
$K$ easily gives an upper bound on $u(K)$, it is rarely clear that this is the best possible. In fact, unknotting
number is probably best understood in a manner that avoids diagrams: it is the minimal number of times
that a knot must pass though itself in a one-parameter family of curves that starts with $K$ and ends with
the unknot.

\noindent {\caps Problem 8.} {\sl Is the unknotting number of a knot algorithmically computable? 
Is there even an algorithm to decide whether a given knot has unknotting number one?}

The first of these questions seems a long way out reach at present. It is more than just a theoretical
question: there exist many explicit knots with unknown unknotting number (see [41] for example).
As a result, many intriguing and clever methods have been developed that provide lower bounds
on unknotting number [38, 39, 42, 60, 74, 76, 79, 94, 96]. These utilise a wide variety of techniques, including the Alexander module [74], the linking form on the double branched cover [60, 94], 4-dimensional gauge theory [96] 
and more recently, Heegaard Floer homology [76, 79].
A major piece of research in this area was the solution by Kronheimer and Mrowka [44] of
Milnor's conjecture: {\sl the unknotting number of the $(p,q)$ torus knot is
$(p-1)(q-1)/2$.}

However, with each new lower bound on unknotting number, there
remain knots with unknotting number that cannot be determined. This will continue
to be the case until Problem 8 is resolved.

A useful way of understanding and approaching this problem is via the use of crossing circles.
Given some crossing in a diagram of a knot $K$, the associated {\sl crossing circle} is a simple closed curve
in the complement of $K$ that encircles the crossing, as shown in Figure 9. There are two crossing
circles associated with each crossing (that are obtained from one another by rotating the diagram
through 90 degrees). One usually focuses on the crossing circle that has zero linking number with $K$.
This bounds a disc which intersects $K$ in two points of opposite sign. 

\vskip 18pt
\centerline{
\includegraphics[width=2.5in]{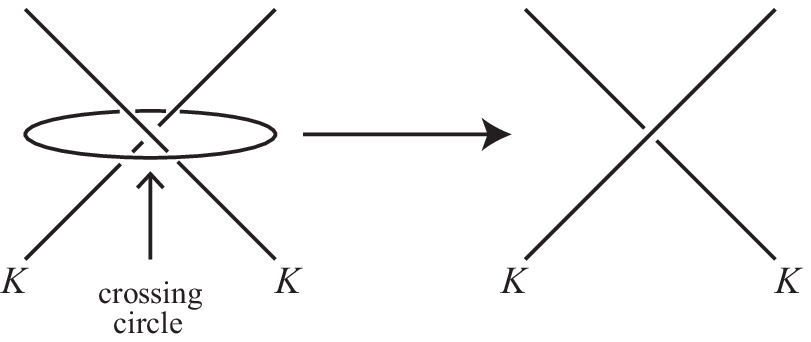}
}
\vskip 6pt
\centerline{Figure 9: A crossing change}

Crossing circles are useful for several
reasons. Firstly, the crossing change is achieved by $\pm 1$ surgery along the crossing
circle, and so one can draw on the well-developed theory of Dehn surgery.
Secondly, it also permits crossings in different diagrams to be compared:
the crossings are {\sl equivalent} if their associated crossing circles
are ambient isotopic in the complement of $K$. Clearly, changing equivalent crossings
results in equivalent knots, and so this is a helpful concept. Some interesting
problems can be phrased in terms of crossing circles.

Although it seems unlikely that we will have an algorithm to determine unknotting number in the near future,
the question of whether unknotting number one knots can be detected is a tantalising one.
Here, the situation seems considerably more promising. To answer this question, it would
suffice to produce a finite list of potential crossing circles with the property that
if $K$ has unknotting number 1, then changing the crossing at one of these
crossing circles would unknot the knot. For one could then check each
of these crossing changes systematically to determine whether they
did unknot the knot. (This would use one of the known solutions to
the unknot recognition problem.) Thus, the following arises naturally.

\noindent {\caps Problem 9.} {\sl If a knot has unknotting number one, are there only
finitely many crossing circles, up to ambient isotopy, that can be used to unknot
the knot?}

Much of what is known about knots with unknotting number one comes from
sutured manifold theory. (See Section 3.5 for a very brief summary of some of this theory.)
It is possible to show that if a knot $K$ has unknotting number 1, 
then its exterior admits a taut sutured manifold hierarchy, {\sl where the penultimate
manifold in this decomposition is a solid torus neighbourhood of the crossing circle plus
possibly some 3-balls.}
In particular, this crossing circle is disjoint from some minimal genus Seifert surface $S$
for $K$. In fact, by pushing this argument further, Scharlemann and Thompson [89]
were able to show that this crossing circle can be isotoped to lie in $S$. This suggests
a route to determining whether a knot has unknotting number 1: find all minimal genus 
Seifert surfaces for $K$, and then search on each such
surface for all possible crossing circles which unknot the knot. The first of these steps
is possible, using normal surface theory, as long as $K$ is not a satellite knot.
However, there are infinitely many isotopy classes of curves on a surface with positive genus, and it is hard
to see how to reduce this to a finite list of crossing circles on which to focus.

\vskip 6pt
\noindent {\caps 6.2. The unlinking and splitting numbers of a link}
\vskip 6pt

There is a natural analogue of unknotting number for links. The {\sl unlinking number} $u(L)$
of a link $L$ is the minimal number of crossing changes required to turn it into the unlink.
Many of the techniques that can be used to analyse unknotting number apply equally
well to unlinking number. But they also often work in the following more general setting.
The {\sl splitting number} $s(L)$ of a link $L$ is the minimal number of crossing changes
required to turn it into a split link.

Unlike the case of unknotting number one, links with splitting number one or unlinking number one are
known to be detectable, under fairly mild hypotheses. It is a recent theorem of the author [57] that
{\sl there is an algorithm to determine whether a hyperbolic link of at least 3 components
has splitting number one}.

This is proved as follows. Suppose that $C$ is a crossing circle, such that
changing this crossing turns a link $L$ into a split link. Just as in Section 6.1,
there is a sutured manifold hierarchy for the exterior of $L$ where the penultimate manifold is a solid
torus neighbourhood of $C$, plus some balls. The key part of the proof is to place this hierarchy
into `normal' form with respect to some initial triangulation of the exterior of $L$ as in [50].
Thus, one obtains a finite computable list of possibilities for $C$. In particular,
the version of Problem 9 in this setting is answered affirmatively. Moreover,
each of these crossing changes can be analysed to determine whether they result in
the unlink. So, under the above hypotheses, {\sl there is also an algorithm to determine whether
the link has unlinking number one}.

\vskip 6pt
\noindent {\caps 6.3. Additivity of unknotting number}
\vskip 6pt

One can also ask how unknotting number behaves with respect to simple knot-theoretic
operations. For example, the following is very far from resolution.

\noindent {\caps Problem 10.} {\sl Does $u(K_1 \sharp K_2) = u(K_1) + u(K_2)$
always hold? Less ambitiously, what about the inequality $u(K_1 \sharp K_2) \geq \max \{ u(K_1), u(K_2) \}$?}

Scharlemann [87] proved that composite knots
have unknotting number at least $2$. But we cannot rule out the following bizarre 
possibility: two knots might each have unknotting number $100$ but their connected sum
has unknotting number $2$!

\vskip 6pt
\noindent {\caps 6.4. Cosmetic crossing changes}
\vskip 6pt

A crossing circle is {\sl nugatory} if it bounds a disc in the complement of
the knot. Clearly changing a nugatory crossing does not alter the knot. Is
this the only way?

\noindent {\caps Problem 11.} {\sl Suppose that a crossing change to a knot $K$ does not alter its
knot type. Let $C$ be the associated crossing circle with zero linking number.
Must $C$ be nugatory?}

This is known when the knot $K$ is the unknot, by work of Scharlemann and
Thompson [89], which relied on Gabai's sutured manifold theory [21]. However, the
general case seems to remain just out of reach.

\vskip 12pt
\centerline{\caps 7. Special classes of knots}
\vskip 6pt

Knots are usually specified by means of a diagram. Typically, this diagram reveals little immediate
information about the knot. For example, a superficially complicated diagram may actually
represent the unknot. However, there are certain types of diagram which force a knot to
have much more structure. It is these that we discuss in this section.

\vskip 6pt
\noindent {\caps 7.1. Alternating knots}
\vskip 6pt

A diagram for a link is {\sl alternating} if, as one travels along each component of the link,
one meets `over' and `under' crossings alternately. A link is {\sl alternating} if
it has an alternating diagram.

\vskip 18pt
\centerline{
\includegraphics[width=1.7in]{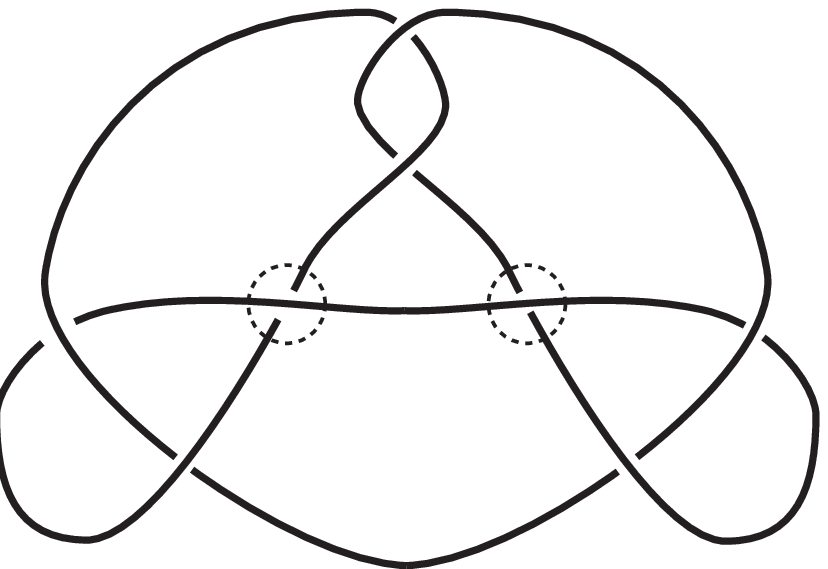}
}
\vskip 6pt
\centerline{Figure 10: A non-alternating diagram (highlighted are successive `over' crossings)}

Some alternating diagrams can be immediately simplified as follows. We say that a 
diagram is {\sl reducible} if, at some crossing, two of the incident regions are actually 
the same region of the diagram. When a diagram is reducible, there is an obvious modification that
can be made to it, which reduces the crossing number by 1.
When the original diagram is alternating, so too is the new diagram. So, when dealing
with alternating links, there is no loss of generality in assuming that its alternating
diagram is not reducible. It is then known as {\sl reduced}.

A great deal of information about a link is evident from an alternating diagram.
Possibly the most striking illustration of this is the following result of
Kauffman [40], Murasugi [72], and Thistlethwaite [97]: {\sl If $D$ is a reduced
alternating diagram of a link $L$, then $c(D) = c(L)$.}

This was proved by means of the Jones polynomial. This link invariant is defined
in terms of some diagram. It is therefore unsurprising that the `complexity' of the
polynomial should be bounded above by the `complexity' of some diagram. In fact,
the {\sl width} of the polynomial, which is the difference in degree between the
highest order term and lowest order term, is always at most the crossing number of any diagram.
But, Kauffman, Murasugi and Thistlethwaite proved that this inequality is an equality
for a reduced alternating diagram $D$. This implies the above result,
because the width of the Jones polynomial is then the crossing number of $D$, and 
the link cannot have a diagram with fewer crossings, because this would force
the width of the polynomial to be at most this number.

A consequence of this theorem is that there is a very rapid and simple algorithm to
decide whether an alternating diagram represents the unknot: {\sl an alternating diagram of the
unknot either has no crossings or is reducible}. In fact, this was already
known, due to earlier techniques of Menasco [68], which relied on a very precise analysis
of surfaces in the exterior of an alternating link. This works particularly
well when the surface has non-negative Euler characteristic. So, by examining spheres
in the complement of the link, one can also detect whether an alternating link is split.
Menasco showed that {\sl an alternating link is split if and only if its alternating
diagram is disconnected}. Similarly, {\sl an alternating knot is prime if and only if
a reduced alternating diagram is prime}. 
Here, a connected diagram is {\sl prime} if and only if any simple closed in the plane
which intersects the link projection transversely in two points away from the crossings divides the plane
into two components, one of which contains no crossings. 

In Menasco's proof, he constructs, using the diagram, a special position for the link
in ${\Bbb R}^3$. The diagram is viewed as lying in the horizontal plane ${\Bbb R}^2 \times \{ 0 \}$.
The link mostly lies in this plane, except near each crossing, where two sub-arcs of the link
take detours above and below the plane of the diagram. More precisely, one places
a small 3-ball at each crossing with centre lying in the plane of the diagram,
and the two arcs involved in the crossing follow the boundary of this ball. (See Figure 11.)
Let ${\Bbb R}^3_+$ denote the set of
points lying above the union of the plane and the balls, and define ${\Bbb R}^3_-$ similarly.
 Menasco showed that
a closed incompressible surface in the link exterior may be arranged so that it
intersects each of these balls in a collection of parallel `saddles', and intersects
${\Bbb R}^3_+$ and ${\Bbb R}^3_-$ in a collection of discs. The boundary of these
discs in ${\Bbb R}^3_+$ is a collection of simple closed curves, which can be viewed as lying within the diagram.
A careful analysis of these curves (particularly, an innermost one in the plane),
together with the condition that the diagram is connected and alternating, 
implies that the surface $S$ has a {\sl meridional compression disc}.
This is a disc $D$ in ${\Bbb R}^3$ such that $D \cap S = \partial D$ and $D \cap K$ is a single
point in the interior of $D$. The existence of such a disc is impossible for a splitting sphere,
for example, and this can be used to establish Menasco's theorem about split
alternating links.

\vskip 18pt
\centerline{
\includegraphics[width=3in]{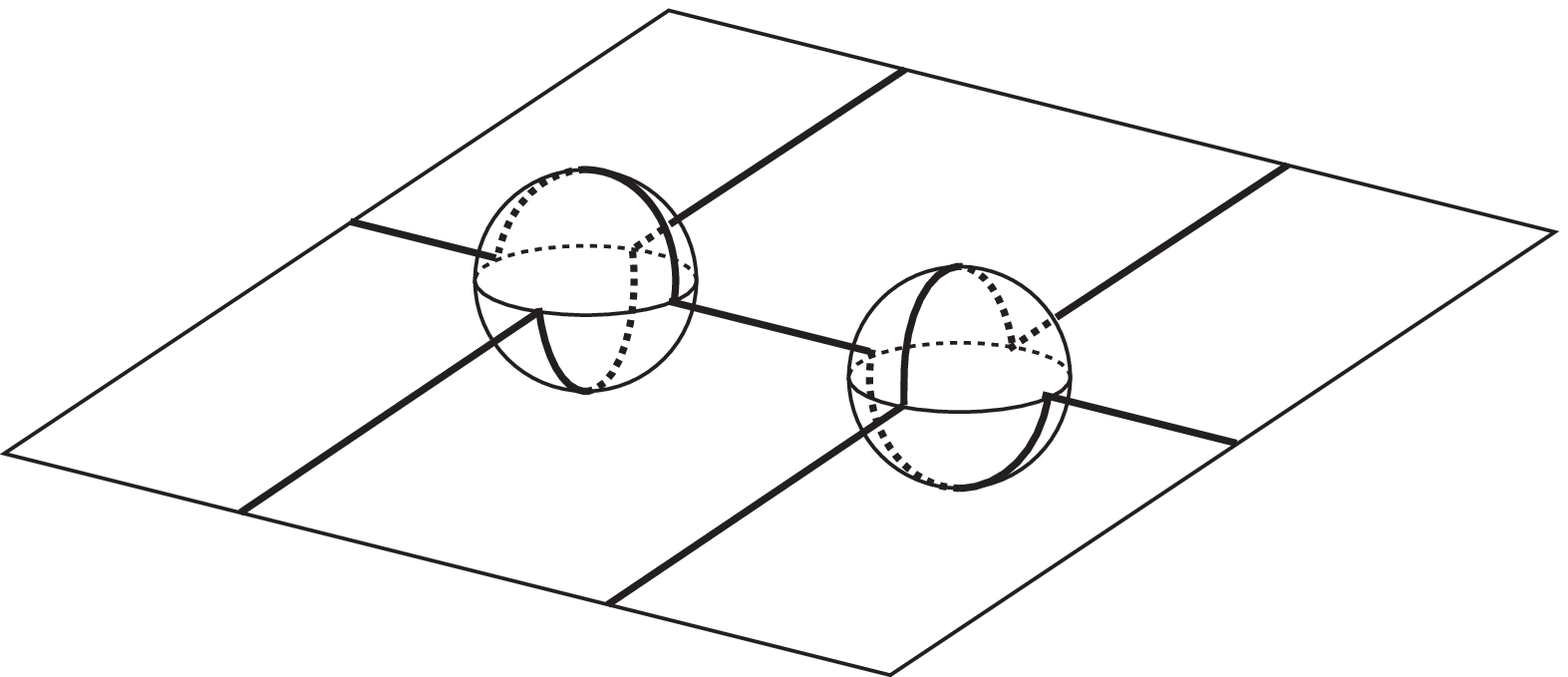}
}
\vskip 6pt
\centerline{Figure 11: Balls near crossings}

A significantly deeper version of this analysis was used to prove the following theorem
of Menasco and Thistlethwaite [69]: {\sl any two reduced alternating diagrams of a link 
differ by a sequence of flypes} (as shown in Figure 12.) This was Tait's flyping conjecture.
The proof of this is very complex,
and so we can only give a very brief overview. Each alternating diagram determines
two checkerboard surfaces. In Menasco and Thistlethwaite's proof, they
examine the surfaces arising from one diagram and how they interact with other
diagram. A careful analysis, together with a surprising intervention of the Jones
polynomial at one point, establishes the existence of a flype which, in some sense,
simplifies these surfaces, and so the theorem is proved by induction.

\vskip 18pt
\centerline{
\includegraphics[width=3in]{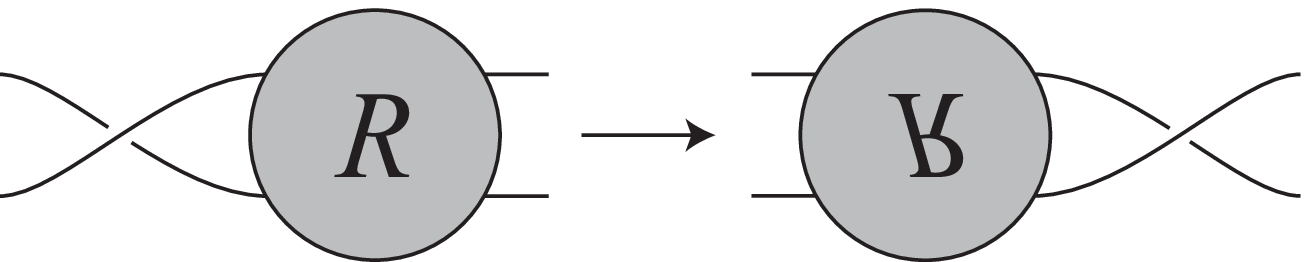}
}
\vskip 6pt
\centerline{Figure 12: Flype}

Although alternating links are now very well-understood, there are some invariants
and properties of these links that are worthy of further analysis. Some of these
are given in the following problem.

\vfill\eject
\noindent {\caps Problem 12.} {\sl Can one detect (or possibly just estimate) the following invariants
of a knot, given an alternating diagram: its bridge number; its tunnel number; its unknotting number?}\footnote{$^3$}{\ninerm 
An algorithm to determine whether an alternating knot has unknotting number one has recently been given by
McCoy [67].}

One may also ask for estimates for more sophisticated invariants. Indeed, in [51], the author has
given a way of reading off the hyperbolic volume of an alternating link's complement, up to a bounded
multiplicative error. While such problems are interesting, they falls outside the remit of this survey, which is focused on much more elementary questions.

\vskip 18pt
\noindent {\caps 7.2. Positive knots and positive braids}
\vskip 6pt

Given the manifest success in understanding the topological properties of alternating knots,
it is natural to consider other possible classes of knot diagrams. A reasonable
way of assessing whether it is a good class of diagrams is whether it satisfies
the minimal requirement that a diagram specifies the unknot if and only if
it `obviously' does so. Another class of diagrams which satisfy this requirement
is positive diagrams.

Positivity is defined for diagrams of oriented links. Positive and negative crossings
are shown in Figure 13. A diagram is said to be {\sl positive} if all its crossings are
positive.

\vskip 18pt
\centerline{
\includegraphics[width=2.3in]{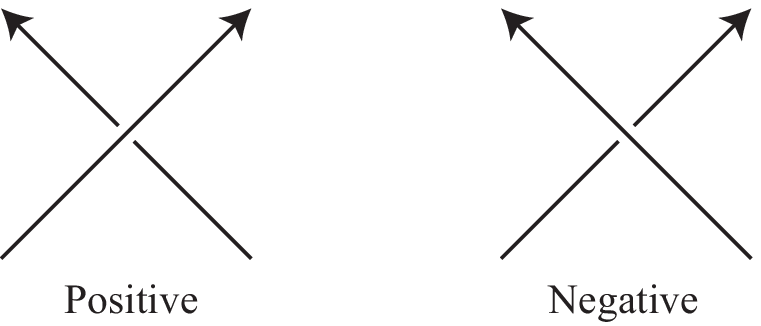}
}
\vskip 6pt
\centerline{Figure 13: Positive and negative crossings}

Cromwell [11] showed that a reduced positive diagram never represents the unknot.
In fact, one can read off the genus of a positive knot from its positive diagram.
This result was significantly generalised by Kronheimer and Mrowka [44],
who showed that the smooth 4-ball genus of a positive knot can be read off
from its positive diagram. Indeed, it is equal to the knot's genus.

However, positive knots are not totally well-behaved. For example,
Stoimenow [93] found a positive knot which has no positive diagram with
minimal crossing number. So, it seems unlikely that one can read
off the crossing number of a positive knot in the way that one can for
alternating knots. Nevertheless, they form an interesting class that is worthy
of further investigation. In particular, it would be instructive to determine 
other topological and/or elementary invariants of positive knots. For example, 
is there a simple relationship between two positive diagrams of the same knot?

One can also consider, more specifically, closed positive braids. These
appear to have even nicer properties than general positive knots. For example,
they are fibred. They may be more amenable to analysis.

\vskip 6pt
\noindent {\caps 7.3. Other interesting classes of diagram}
\vskip 6pt

There are several other related classes of diagrams. {\sl Homogeneous}
diagrams form a common generalisation of positive diagrams and alternating
diagrams. Cromwell [11] in fact established his formula for knot genus
not just for positive knots but homogeneous ones.

Adequate diagrams form another interesting class. These are defined by considering
`resolutions' for each crossing, which remove the crossing. Each crossing may be removed in one of two ways,
a `plus' or a `minus' resolution. If one resolves every crossing, the result
is a collection of planar curves. At each crossing, the two new arcs may belong
to the same simple curve or two distinct ones. A diagram is {\sl plus-adequate}
if, when all the crossings are resolved in a `plus' way, each resulting simple closed
curve runs through the remnants of each crossing at most once. {\sl Minus-adequacy} is defined
similarly using the `minus' resolution at each crossing. A diagram is {\sl adequate}
if it is both plus-adequate and minus-adequate. Adequate knots have well-behaved Jones
polynomial. Lickorish and Thistlethwaite [61] using this to prove that {\sl if a knot has
a reduced adequate diagram with at least one crossing, then it is
not the unknot.} But adequacy has more topological applications, as shown in
the recent monograph of Futer, Kalfagianni and Purcell [20]. Adequate links surely merit
further investigation.

\vskip 18pt
\centerline{\caps 8. Epilogue}
\vskip 6pt

The theory of 3-manifolds has seen remarkable progress over the past few years,
with the solution of many key conjectures, for example, the Poincar\'e conjecture [81, 82, 83]
and the Virtually Haken Conjecture [2], to name just two. One might be led to the conclusion
that the related field of knot theory might be devoid of interesting unsolved problems.
The purpose of the present article is to emphasise that this is not the case. There remain
many fundamental and interesting challenges ahead.

\vskip 18pt
\centerline{\caps References}
\vskip 6pt

{\parskip = 0.05truein 
\baselineskip = 0.17truein
{\ninerm

\item{1.} {\smallcaps I. Agol}, {\smallsl Knot genus is NP}, http://www.math.uic.edu/ $\sim$ agol/research.html

\item{2.} {\smallcaps I. Agol}, {\smallsl The virtual Haken conjecture.} (with an appendix by I. Agol, D. Groves and J. Manning)
Doc. Math. 18 (2013) 1045--1087. 

\item{3.} {\smallcaps I. Agol, J. Hass, W. Thurston}, {\smallsl The computational complexity of knot genus and spanning area},
Trans. Amer. Math. Soc. 358 (2006) 3821--3850.

\item{4.} {\smallcaps D. Bennequin,} {\smallsl Entrelacements et \'equations de Pfaff}, 
Ast\'erisque 107--108 (1983) 87--161.

\item{5.} {\smallcaps J. Birman, T. Brendle}, {\smallsl Braids: a survey.} 
Handbook of knot theory, 19--103, Elsevier B. V., Amsterdam, 2005.

\item{6.} {\smallcaps J. Birman, M. Hirsch,} {\smallsl A new algorithm for recognizing the unknot}, 
Geometry and Topology 2 (1998) 178--220.

\item{7.} {\smallcaps J. Birman, W. Menasco}, {\smallsl Studying links via closed braids V: Closed braid
representatives of the unlink,} Trans. Amer. Math. Soc. 329 (1992) 585--606.

\item{8.} {\smallcaps M. Brittenham}, {\smallsl Essential laminations and Haken normal form.}
Pacific J. Math. 168 (1995), no. 2, 217--234.

\item{9.} {\smallcaps A. Coward, M. Lackenby}, {\smallsl An upper bound on Reidemeister moves},
Amer. J. Math. 136 (2014), no. 4, 1023--1066.

\item{10.} {\smallcaps D. Coulson, O. Goodman, C. Hodgson, W. Neumann}, {\smallsl  Snap}, Available at \break
http://www.ms.unimelb.edu.au/$\sim$snap/

\item{11.} {\smallcaps P. Cromwell,} {\smallsl Homogeneous links.} J. London Math. Soc. (2) 39 (1989), no. 3, 535--552.

\item{12.} {\smallcaps P. Cromwell,} {\smallsl Embedding knots and links in an open book I: basic properties}, 
Topology and its Applications 64 (1995) 37--58.

\item{13.} {\smallcaps M. Culler, P. Shalen,} {\smallsl Bounded, separating, incompressible surfaces in knot manifolds,}
Invent. Math. 75 (1984) 537--545.

\item{14.} {\smallcaps F. Dahmani, D. Groves,} {\smallsl The isomorphism problem for toral relatively 
hyperbolic groups,} Publ. Math. Inst. Hautes \'Etudes Sci. No. 107 (2008), 211--290.

\item{15.} {\smallcaps Y. Diao,} {\smallsl The additivity of crossing numbers,} 
Journal of knot theory and its Ramifications 13 (2004) 857--866.

\item{16.} {\smallcaps S. Donaldson,} {\smallsl An Application of Gauge Theory to Four Dimensional Topology}, 
J. Differential Geom. 18 (183) 279--315.

\item{17.} {\smallcaps I. Dynnikov,} {\smallsl Arc-presentations of links: monotonic simplification,}
Fund. Math. 190 (2006), 29--76.

\item{18.} {\smallcaps D. Epstein, R. Penner}, {\smallsl Euclidean decompositions of noncompact hyperbolic manifolds},
 J. Differential Geom. 27 (1988) 67--80.

\item{19.} {\smallcaps A. Floer}, {\smallsl An instanton-invariant for 3-manifolds}, 
Comm. Math. Phys. 118 (2) 215--240.

\item{20.} {\smallcaps D. Futer, E. Kalfagianni, J. Purcell},
{\smallsl Guts of surfaces and the colored Jones polynomial},
Lecture Notes in Mathematics, 2069. Springer, Heidelberg, 2013

\item{21.} {\smallcaps D. Gabai,} {\smallsl Foliations and the topology of $3$-manifolds.}
J. Differential Geom. 18 (1983), no. 3, 445--503.

\item{22.} {\smallcaps D. Gabai,} {\smallsl Essential laminations and Kneser normal form. }
J. Differential Geom. 53 (1999), no. 3, 517--574. 

\item{23.} {\smallcaps C. Gordon, J. Luecke,} {\smallsl Knots are determined by their complements.}
J. Amer. Math. Soc. 2 (1989), no. 2, 371--415. 

\item{24.} {\smallcaps W. Haken,} Theorie der Normalfl\"achen. Acta Math. 105 (1961) 245--375.

\item{25.} {\smallcaps W. Haken,} {\smallsl Some results on surfaces in $3$-manifolds.}
Studies in Modern Topology pp. 39--98 Math. Assoc. Amer. (1968)

\item{26.} {\smallcaps J. Hass, J. Lagarias,} {\smallsl The number of Reidemeister moves needed for unknotting.}
J. Amer. Math. Soc. 14 (2001), no. 2, 399--428 

\item{27.} {\smallcaps J. Hass, J. Lagarias, N. Pippenger,} {\smallsl The computational complexity of knot and 
link problems.} J. ACM 46 (1999), no. 2, 185--211.

\item{28.} {\smallcaps J. Hass, T. Nowik,} {\smallsl Unknot diagrams requiring a quadratic number of 
Reidemeister moves to untangle.} Discrete Comput. Geom. 44 (2010), no. 1, 91--95.

\item{29.} {\smallcaps J. Hass, J. Snoeyink, W. Thurston,} {\smallsl The size of spanning disks for polygonal curves.}
Discrete Comput. Geom. 29 (2003) 1--17.

\item{30.} {\smallcaps G. Hemion}, {\smallsl  On the classification of homeomorphisms of $2$-manifolds 
and the classification of $3$-manifolds.} Acta Math. 142 (1979), no. 1-2, 123--155. 

\item{31.} {\smallcaps J. Hoste, M. Thistlethwaite,} {\smallsl Knotscape},
http://www.math.utk.edu/$\sim$morwen/knotscape.html

\item{32.} {\smallcaps J. Hoste, M. Thistlethwaite, J. Weeks,} 
{\smallsl The first 1,701,936 knots.} Math. Intelligencer 20 (1998), no. 4, 33Ð48.

\item{33.} {\smallcaps W. Jaco, P. Shalen,} {\smallsl Seifert fibered spaces in 3-manifolds.}
Mem. Amer. Math. Soc. 21 (1979), no. 220.

\item{34.} {\smallcaps W. Jaco, J. Tollefson,} {\smallsl Algorithms for the complete decomposition of a 
closed 3-manifold,} Illinois J. Math. 39 (1995), no. 3, 358--406.

\item{35.} {\smallcaps F. Jaeger, D. Vertigana, D. Welsh}, {\smallsl
On the computational complexity of the Jones and Tutte polynomials},
Math. Proc. Cambridge Phil. Soc. 108 (1990) 35--53.

\item{36.} {\smallcaps K. Johannson,} {\smallsl Homotopy equivalences of 3-manifolds with boundaries.}
Lecture Notes in Mathematics, 761. Springer, Berlin, 1979.

\item{37.} {\smallcaps V. Jones}, {\smallsl A polynomial invariant for knots via von Neumann algebras.}
Bull. Amer. Math. Soc. (N.S.) 12 (1985), no. 1, 103--111.

\item{38.} {\smallcaps T. Kawamura}, {\smallsl The unknotting numbers of $10_{139}$ and $10_{152}$ are $4$.}
Osaka J. Math. 35 (1998), no. 3, 539--546.

\item{39.} {\smallcaps T. Kanenobu, H. Murakami,} {\smallsl Two-bridge knots with unknotting number one.}
Proc. Amer. Math. Soc. 98 (1986) 499--502.

\item{40.} {\smallcaps L. Kauffman}, {\smallsl State models and the Jones polynomial.}
Topology 26 (1987) 395--407. 

\item{41.} {\smallcaps R. Kirby}, {\smallsl Problems in low dimensional manifold theory.}
Algebraic and geometric topology (Proc. Sympos. Pure Math., Stanford Univ., Stanford, Calif., 1976), 
Part 2, pp. 273-312, Proc. Sympos. Pure Math., XXXII, Amer. Math. Soc., Providence, R.I., 1978. 

\item{42.} {\smallcaps T. Kobayashi,} {\smallsl Minimal genus Seifert surface for unknotting number 1 knots.} 
Kobe J. Math. 6 (1989) 53--62.

\item{43.} {\smallcaps P. Koiran,} {\smallsl HilbertÕs Nullstellensatz is in the polynomial hierarchy,} 
J. Complexity 12 (1996), no. 4, 273Ð286, DIMACS TR 96-27, Special issue for FOCM 1997.

\item{44.} {\smallcaps P. Kronheimer, T. Mrowka,}
{\smallsl Gauge theory for embedded surfaces. I.}
Topology 32 (1993), no. 4, 773--826. 

\item{45.} {\smallcaps P. Kronheimer, T. Mrowka}, {\smallsl Dehn surgery, the fundamental group and SU(2),}
Math. Res. Lett. 11 (2004) 741--754.

\item{46.} {\smallcaps P. Kronheimer, T. Mrowka}, {\smallsl Monopoles and three-manifolds.} 
New Mathematical Monographs, 10. Cambridge University Press, Cambridge, 2007.

\item{47.} {\smallcaps P. Kronheimer, T. Mrowka}, {\smallsl Khovanov homology is an unknot-detector,}
Publ. Math. Inst. Hautes \'Etudes Sci. 113 (2011) 97--208.

\item{48.} {\smallcaps G. Kuperberg}, {\smallsl Knottedness is in NP, modulo GRH}, 
Adv. Math. 256 (2014), 493--506.

\item{49.} {\smallcaps G. Kuperberg}, {\smallsl Algorithmic homeomorphism of 3-manifolds as a corollary of geometrization}, \hfill\break
arXiv:1508.06720

\item{50.} {\smallcaps M. Lackenby}, {\smallsl Exceptional surgery curves in triangulated 3-manifolds,}
Pacific J. Math. 210 (2003), 101--163.

\item{51.} {\smallcaps M. Lackenby}, {\smallsl The volume of hyperbolic alternating link complements,}
Proc. London Math. Soc. 88 (2004) 204--224.

\item{52.} {\smallcaps M. Lackenby}, {\smallsl  The crossing number of composite knots}, 
J. Topology 2 (2009) 747-768. 

\item{53.} {\smallcaps M. Lackenby}, {\smallsl Core curves of triangulated solid tori}, 
arXiv:1106.2934, Trans. Amer. Math. Soc. 366 (2014), no. 11, 6027--6050.

\item{54.} {\smallcaps M. Lackenby}, {\smallsl  The crossing number of satellite knots}, 
arXiv:1106.3095, Alg. Geom. Top. 14 (2014), no. 4, 2379--2409.

\item{55.} {\smallcaps M. Lackenby}, {\smallsl A polynomial upper bound on Reidemeister moves},
 Ann. of Math. 182 (2015), no. 2, 491--564.

\item{56.} {\smallcaps M. Lackenby}, {\smallsl The efficient certification of knottedness and Thurston norm},
Preprint.

\item{57.} {\smallcaps M. Lackenby}, {\smallsl Links with splitting number one}, In preparation.

\item{58.} {\smallcaps M. Lackenby}, {\smallsl A polynomial upper bound on Reidemeister moves for each link type},
In preparation.

\item{59.} {\smallcaps J. Lagarias, A. Odlyzko,} {\smallsl Effective versions
of the Chebotarev density theorem,} Algebraic number Þelds: L-functions and Galois properties 
(Proc. Sympos., Univ. Durham, 1975), Academic Press, 1977, pp. 409--464.

\item{60.} {\smallcaps W. B. R. Lickorish,} {\smallsl The unknotting number of a classical knot.} Contemp.
Math. 44 (1985) 117--121.

\item{61.} {\smallcaps W. B. R. Lickorish, M. Thistlethwaite}, {\smallsl Some links with nontrivial polynomials 
and their crossing-numbers}, Comment. Math. Helv. 63 (1988), no. 4, 527--539. 

\item{62.} {\smallcaps A. Mal'cev,} {\smallsl On isomorphic matrix representations of infinite groups,} 
Mat. Sb. 8 (50) (1940), 405--422.

\item{63.} {\smallcaps J. Manning}, {\smallsl Algorithmic detection and description of hyperbolic structures 
on closed 3-manifolds with solvable word problem,} Geom. Topol. 6 (2002), 1--25.

\item{64.} {\smallcaps C. Manolescu, P. Ozsv\'ath, S. Sarkar,} {\smallsl A combinatorial description of knot Floer homology.}
Ann. of Math. (2) 169 (2009), no. 2, 633--660.

\item{65.} {\smallcaps A. Markov}, {\smallsl The unsolvability of the homeomorphy problem},
Dokl. Akad. Nauk SSSR, 121 (1958) 218--220.

\item{66.} {\smallcaps S. Matveev}, {\smallsl Algorithmic topology and classification of 3-manifold},
Algorithms and Computation in Mathematics, 9. Springer, Berlin, 2007.

\item{67.} {\smallcaps D. McCoy,} {\smallsl Alternating knots with unknotting number one},
arXiv:1312.1278

\item{68.} {\smallcaps W. Menasco}, {\smallsl Closed incompressible surfaces in alternating knot and link complements.}
Topology 23 (1984) 37--44.

\item{69.} {\smallcaps W. Menasco, M. Thistlethwaite}, {\smallsl The classification of alternating links.}
Ann. of Math. (2) 138 (1993) 113--171.

\item{70.} {\smallcaps A. Mijatovic}, {\smallsl Simplifying triangulations of $S^3$.} 
Pacific J. Math. 208 (2003), no. 2, 291--324.

\item{71.} {\smallcaps A. Mijatovic}, {\smallsl Simplical structures of knot complements,} 
Math. Res. Lett. 12 (2005) 843--856.

\item{72.} {\smallcaps K. Murasugi}, {\smallsl Jones polynomials and classical conjectures in knot theory.}
Topology 26 (1987), no. 2, 187Ð194. 

\item{73.} {\smallcaps K. Murasugi,} {\smallsl On the braid index of alternating links.}
Trans. Amer. Math. Soc. 326 (1991) 237--260. 

\item{74.} {\smallcaps  Y. Nakanishi,} {\smallsl A note on unknotting number.} 
Math. Sem. Notes Kobe Univ. 9 (1981) 99--108.

\item{75.} {\smallcaps W. Neumann, G. Swarup,} {\smallsl Canonical decompositions of 3-manifolds.}
Geom. Topol. 1 (1997), 21--40.

\item{76.} {\smallcaps B. Owens,} {\smallsl Unknotting information from Heegaard Floer homology.}
Adv. Math. 217 (2008), no. 5, 2353--2376.

\item{77.} {\smallcaps P. Ozsv\'ath and Z. Szab\'o},
{\smallsl Holomorphic disks and topological invariants for closed three-manifolds.}
Ann. of Math. (2) 159 (2004), no. 3, 1027--1158. 

\item{78.} {\smallcaps P. Ozsv\'ath and Z. Szab\'o}, {\smallsl
Holomorphic disks and genus bounds},
Geom. Topol. 8 (2004) 311--334.

\item{79.} {\smallcaps P. Ozsv\'ath and Z. Szab\'o},
{\smallsl Knots with unknotting number one and Heegaard Floer homology}, 
Topology 44 (2005), no. 4, 705--745.

\item{80.} {\smallcaps U. Pachner} {\smallsl 
P.L. homeomorphic manifolds are equivalent by elementary shellings,} 
European Journal of Combinatorics 12 (1991) 129--145.

\item{81.} {\smallcaps G. Perelman,} {\smallsl The entropy formula for the 
Ricci flow and its geometric applications,} Preprint,
arxiv:math.DG/0211159

\item{82.} {\smallcaps G. Perelman,} {\smallsl  Ricci flow with surgery on three-manifolds,}
Preprint, arxiv:math.DG/0303109

\item{83.} {\smallcaps G. Perelman,} {\smallsl Finite extinction time for the 
solutions to the Ricci flow on certain three-manifolds,} Preprint,
arxiv:math.DG/0307245 

\vfill\eject
\item{84.} {\smallcaps V. Pratt.} {\smallsl Every prime has a succinct certificate.} SIAM Journal on Computing, vol.4, pp.214--220. 1975.

\item{85.} {\smallcaps K. Reidemeister,} {\smallsl Knotten und Gruppen,} Abh. Math. Sem. Univ. Hamburg 5 (1927), 7--23.

\item{86.} {\smallcaps S. Sarkar, J. Wang,} {\smallsl An algorithm for computing some Heegaard Floer homologies.}
Ann. of Math. (2) 171 (2010), no. 2, 1213--1236.

\item{87.} {\smallcaps M. Scharlemann}, {\smallsl Unknotting number one knots are prime.}
Invent. Math. 82 (1985), no. 1, 37--55.

\item{88.} {\smallcaps M. Scharlemann,} {\smallsl Sutured manifolds and generalized Thurston norms.}
J. Differential Geom. 29 (1989), no. 3, 557--614.

\item{89.} {\smallcaps M. Scharlemann, A. Thompson,} {\smallsl Link genus and the Conway moves.}
Comment. Math. Helv. 64 (1989), no. 4, 527--535

\item{90.} {\smallcaps M. Scharlemann}, {\smallsl Crossing changes. Knot theory and its applications. }
Chaos Solitons Fractals 9 (1998) 693--704.

\item{91.} {\smallcaps P. Scott, H. Short,} {\smallsl The homeomorphism problem for closed 3-manifolds.}
Algebr. Geom. Topol. 14 (2014) 2431--2444. 

\item{92.} {\smallcaps Z. Sela,} {\smallsl The isomorphism problem for hyperbolic groups. I.}
Ann. of Math. (2) 141 (1995) 217--283. 

\item{93.} {\smallcaps A. Stoimenow}, {\smallsl On the crossing number of positive knots and braids 
and braid index criteria of Jones and Morton-Williams-Franks}, Trans. Amer. Math. Soc. 354 (2002), 
3927--3954. 

\item{94.} {\smallcaps A. Stoimenow,} {\smallsl Polynomial values, the linking form and unknotting numbers.}
Math. Res. Lett. 11 (2004), no. 5-6, 755--769.

\item{95.} {\smallcaps C. Sundberg, M. Thistlethwaite,} {\smallsl 
The rate of growth of the number of prime alternating links and tangles}, Pacific J. Math. 
182 (1998) 329--358.

\item{96.} {\smallcaps T. Tanaka,} {\smallsl Unknotting numbers of quasipositive knots},
Topology Appl. 88 (1998), no. 3, 239--246.

\item{97.} {\smallcaps M. Thistlethwaite}, {\smallsl A spanning tree expansion of the Jones polynomial.}
Topology 26 (1987)  297--309. 

\item{98.} {\smallcaps M. Thistlethwaite}, {\smallsl On the structure and scarcity of alternating links and tangles.}
J. Knot Theory RamiÞcations 7 (1998) 981--1004.

\item{99.} {\smallcaps A. Thompson}, {\smallsl Thin position and the recognition problem for $S^3$.}
Math. Res. Lett. 1 (1994), no. 5, 613--630. 

\item{100.} {\smallcaps W. Thurston}, {\smallsl Three-dimensional manifolds, Kleinian groups 
and hyperbolic geometry.} Bull. Amer. Math. Soc. (N.S.) 6 (1982), no. 3, 357--381.

\item{101.} {\smallcaps W. Thurston,} {\smallsl A norm for the homology of 3-manifolds.}
Mem. Amer. Math. Soc. 59 (1986), no. 339, i--vi and 99--130. 

\item{102.} {\smallcaps J. Weeks}, {\smallsl Snappea}, http://www.geometrygames.org/SnapPea/

\item{103.} {\smallcaps J. Weeks,} {\smallsl 
Convex hulls and isometries of cusped hyperbolic 3-manifolds.} Topology Appl. 52 (1993), no. 2, 127--149.

\item{104.} {\smallcaps P. Weinberger}, {\smallsl Finding the number of factors of a polynomial,} 
J. Algorithms 5 (1984) 180--186.

\item{105.} {\smallcaps E. Witten}, {\smallsl Quantum Field  Theory  and the  Jones Polynomial},
Commun. Math.  Phys.  121 (1989) 351--399.

\item{106.} {\smallcaps E. Witten,} {\smallsl Monopoles and four-manifolds}, 
Math. Res. Lett. 1 (6) 769--796.

}}

\vskip 12pt
\+ Mathematical Institute, University of Oxford, \cr
\+ Radcliffe Observatory Quarter, Woodstock Road, Oxford OX2 6GG, United Kingdom. \cr

\end